\documentclass[12pt]{elsarticle}
\usepackage{cancel}
\usepackage{caption}
\usepackage{color}
\usepackage{lscape}
\usepackage{afterpage}
\usepackage{pstricks}
\usepackage{pst-plot}
\usepackage{longtable}
\usepackage{dcolumn}
\usepackage{pst-node}
\usepackage{amsmath}
\usepackage{amssymb}
\usepackage{amsthm}
\usepackage{multirow}
\usepackage{colortab}
\usepackage{color}
\usepackage{rotating}
\usepackage{array}
\usepackage{slashbox}
\usepackage{colortbl}
\usepackage{subfigure}
\usepackage{textcomp}
\usepackage{pstricks, pst-node}
\usepackage{pst-all}
\usepackage{algorithm,algorithmic}
\usepackage{url}
\usepackage{soul}
\usepackage{hyperref}
\psset{arrows=->, labelsep=3pt, mnode=circle}
\usepackage{lineno}

\usepackage{tabularx}

\usepackage{eurosym}

\newfont{\rams}{msbm10 scaled\magstep1}
\newcommand{\rea}{\mbox{\rams \symbol{'122}}}

\newenvironment{resumeT}{\begin{list}{}{\setlength{\rightmargin}{\leftmargin}}\item[]
{\centering {\bf \it~~~}
\par}\item[]\ignorespaces}{\unskip\end{list}}

\begin{document}
\title{Multiobjective Combinatorial Optimization with Interactive Evolutionary Algorithms: the case of facility location problems}

\author[Ven,Por]{\rm Maria Barbati}
\ead{maria.barbati@unive.it}
\author[Eco]{\rm Salvatore Corrente}
\ead{salvatore.corrente@unict.it}
\author[Eco,Por]{\rm Salvatore Greco}
\ead{salgreco@unict.it}

\address[Eco]{Department of Economics and Business, University of Catania, Corso Italia, 55, 95129  Catania, Italy}
\address[Ven]{Universit\'{a} C\'{a} Foscari Venezia, Department of Economics, 30121 Venice, Italy}
\address[Por]{Portsmouth Business School, Centre of Operations Research and Logistics (CORL), University of Portsmouth, Portsmouth, United Kingdom}

\date{}
\maketitle

\vspace{-1cm}

\begin{resumeT}
{\large {\bf Abstract:}}  
We consider multiobjective combinatorial optimization problems handled by means of preference driven efficient heuristics. They look for the most preferred part of the Pareto front on the basis of some preferences expressed by the Decision Maker during the process. In general, what is searched for in this case is the Pareto set of efficient solutions. This is a problem much more difficult than optimizing a single objective function. Moreover, obtaining the Pareto set does not mean that the decision problem is solved since one or some of the solutions have to be chosen. Indeed, to make a decision, it is necessary to determine the most preferred solution in the Pareto set, so that it is also necessary to elicit the preferences of the user. In this perspective, what we are proposing  can be seen as the first structured methodology in facility location problems to search optimal solutions taking into account preferences of the user. With this aim, we approach facility location problems using a recently proposed interactive evolutionary multiobjective optimization procedure called NEMO-II-Ch. NEMO-II-Ch is applied to a real world multiobjective location problem with many users and many facilities to be located. Several simulations considering different fictitious users have been performed. The results obtained by NEMO-II-Ch are compared with those got by three algorithms which know the user's true value function that is, instead, unknown to NEMO-II-Ch. They show that in many cases NEMO-II-Ch finds the best subset of locations more quickly than the methods knowing, exactly, the whole user's true preferences.


\vspace{0,3cm}
\noindent{\bf Keywords}: {Multiobjective Optimization, Combinatorial Optimization, Preferences, NEMO, Facility Location problems}
\end{resumeT}

\pagenumbering{arabic}

\section{Introduction}\label{Intro}
Multiple Objective Combinatorial Optimization (MOCO) problems (for a survey see \cite{Ehrgott:2000}) are very complex and difficult to be solved. They can be approached with different aims but, in general, one focuses on the computation of all the efficient solutions (see \cite{Serafini1987} for a discussion on the different concepts of solutions of  a MOCO problem). In general, the number of efficient solutions grows exponentially with the size of the problem \cite{Ehrgott:2000,teghem2000interactive}. This, together with the intrinsic complexity related to the ``nonsmothness" of the optimization problems, requires a huge computational effort, much greater than that one involved in the resolution of the single objective cases \cite{alves2007review}. The high number of efficient solutions and the required very high computational effort are considered the main bottleneck of the MOCO problem \cite{alves2007review,coutinho2012bi}. These considerations have triggered the development of a certain number of approaches using heuristics that are able to determine an approximation of the whole set of nondominated or efficient solutions involving less computational effort than that one involved in exact algorithms \cite{Ehrgott:2008}. However, observe that if from a mere theoretical point of view these can be seen as the main critical issues of a MOCO problem, from the point of view of real life applications, there are other difficulties. Indeed, it would be hard to say that a problem has been solved even in case the whole set of efficient solutions has been computed. This set can  contain even several thousands of elements so that, finally, the Decision Maker (DM) who should choose one or some of them could feel himself lost \cite{alves2007review,coutinho2012bi}. Therefore, beyond the technical limitations related to the computational aspects, there are a little more practical questions related to the support given to the DM. From this point of view, the algorithms can take advantage from the integration of preferences expressed by the DM guiding the search to the part of the Pareto front most interesting for him.  {Considering different moments in which the DM is asked to provide his preferences, in literature one distinguishes between a priori, interactive and a posteriori methods} \cite{Ehrgott:2000}:

\begin{itemize}
\item in \textit{a priori} methods, the preferences of the DM are articulated at the beginning of the process,
\item in the \textit{interactive}  methods, the DM expresses his preferences during the search,
\item in the \textit{a posteriori} methods the DM is presented with the set of all efficient solutions that is therefore analyzed w.r.t. his preferences.
\end{itemize}

On the one hand, the use of a priori methods asks the DM to define at the beginning of the procedure his preferences that are translated by some particular utility function. This assumes that the DM is rational and that his decisions are taken on the basis of some pre-existing preferences that, consequently, have only to be discovered. However, this is not always true since the DM not only is uncertain about his preferences at the beginning of the search but, even more, these a priori preferences are in general absent and have to be constructed during the decision process \cite{roy1987meaning,roy1993decision}. \\
On the other hand, in the a posteriori methods the DM is often presented with many solutions. This approach has some drawbacks too since:
\begin{itemize}
\item the DM has to choose the best solution(s) analyzing the tradeoffs among objectives \cite{coutinho2012bi},
\item showing the whole set of solutions can cause an information overload on the DM who may have difficulty in selecting the best one(s) \cite{karasakal2008interactive}.
\end{itemize}

\noindent From what we said above, using interactive methods seems the best choice \cite{CoCo02,Serafini1987,VanVeldhuizenPhD}. Therefore, for MOCO problems, a reasonable approach seems the use of specific heuristics that instead of approximating the whole set of efficient solutions, look for some efficient solutions being the most preferred by the DM. This implies that the heuristics used to explore the feasible set of solutions incorporate few preference information supplied by the DM permitting to drive the search towards some regions of the Pareto front containing the most preferred solutions for the DM. This is possible using  some  recently proposed heuristics \cite{BGSZ-IEEE-TEC-2014} that combine the capacity of a ``smart" exploration of the set of feasible solutions (typical of multiobjective optimization oriented heuristics such as NSGA-II \cite{debieee} or SPEA \cite{zlt2002a}) with the capacity to build a decision model representing DM's preferences (typical of some Multiple Criteria Decision Aiding (MCDA) approaches such as ordinal regression \cite{GrecoEhrgottFigueira2016,jacquet1982assessing}). One example of such composite methodology is given by the recently proposed NEMO-II-Ch algorithm \cite{NEMOIICh} that combines the search procedure of NSGA-II with the preference representation obtained by the nonadditive robust ordinal regression \cite{angilella2010non}. This approach that drives the search of optimal solutions guided by a preference model incorporating the preferences expressed by the DM, seems us a very promising approach to MOCO problems in real life applications. Indeed, it can give appropriate answers to all the limitations of MOCO problems that we have described:
\begin{itemize}
\item	it handles the big number of efficient solutions of a MOCO problem by looking only to small subsets of efficient solutions that are well appreciated by the DM,
\item it handles the computational effort by using heuristics that have proved to be very effective in complex multiobjective problems,
\item it handles the request of a decision support by driving the whole search algorithm by the preferences step by step expressed by the DM in an interactive procedure.
\end{itemize}

To test the usefulness of such an approach in this paper we consider a typical MOCO problem that is the Facility Location Problem (FLP) \cite{Farahanisurevey:2010}.

In FLPs we aim to locate a set of facilities in a space optimizing some objective functions and satisfying some constraints. Historically FLPs have been modeled using a mono objective approach in which a single objective function has been adopted. Many have been the contributions in this sense with a multitude of  objectives adopted \cite{knuth:art} for describing several and very different applications \cite{LocationScience}. However, in the real world, DMs deal with several conflicting objectives at the same time so that it is advisable that also the adopted algorithms take into account a multiobjective formulation of the problem at hand \cite{Ehrgott:2000,Roy_book_greco_2016}.

The classical approach for choosing the position of a facility consists in describing a function of the distances between the potential users of the facility and the facility itself \cite{eiseltLaporte:1995}. The objective function becomes a linear mathematical expression of the distances to optimize. Introducing several constraints, a combinatorial optimization model is built and the optimal solution can be found from the resolution of the model, as described in what is considered the seminal paper by \cite{Hakimi:1964}. Therefore the main aim becomes the theoretical development and the description of properties of the models  and their solutions \cite {LocationScience}. Some reviews gather the basic knowledge on location science as in  \cite{Daskin}  and in the recent books of \cite{Drezner:2004a} and \cite{LocationScience}.

Multiple Objectives Facility Location Problems (MOFLPs) have captured attention from the researchers especially in the last decade. Many objectives can be used: from the classical distance related objectives, to the environmental and ecological criteria (for a list see \cite{Farahanisurevey:2010}). The majority of the methodologies aims to find the whole Pareto front or a part of it implying a considerable computational effort \cite{alves2007review}. To this aim several methodologies can be adopted: from exact approaches (e.g. \cite{Hamacher:2002,Ohsawa:2003})  to  multiobjective evolutionary algorithms (e.g.\cite{deb-book-01}) for complex problems.

Behind all these approaches there is the strong assumption that the DM is able to select the alternative that is the best for him which implies that the DM has clear and well defined  preferences and is completely rational. In most of the practical problems these assumptions are not very realistic \cite{alves2007review,miettinen2008introduction}. Moreover, very few papers take into account directly the opinion of the DMs. Often the objectives described are derived from considerations related to the specific problem without investigating further the opinion of DMs. Therefore, handling complex MOFLPs with an optimization algorithm guided by DM's preferences seems an interesting approach to be explored and, in this perspective, our contribution can be considered the first structured methodology in MOFLPs to search optimal solutions taking into account preferences of the user. For this reason we propose to deal with MOFLPs by using NEMO-II-Ch. In this way, interactively the DM provides his preferences on some pairs of possible facility locations assignments guiding therefore the search to the part of the Pareto front most interesting w.r.t his preferences and avoiding to loose time in looking for other solutions not matching his expectation. We shall present our proposal considering a case study introduced in \cite{drezner2006multi}.

To underline the efficiency of NEMO-II-Ch to MOCO problems, we simulated different user's value functions and we compared the algorithm's performances with the ones of three algorithms, denoted by EA-UVF \citep{NEMOIICh}, EA-UVF1 and EA-UVF2, based on the knowledge of the user's true value function. We observed that, quite often, NEMO-II-Ch performs better than the algorithms knowing the user's preferences. To test how the convergence of NEMO-II-Ch to the preferred solution is dependent on the preference information provided by the DM with the related cognitive burden, we considered three different variants asking him to compare one pair of solutions every 5, 10 and 20 generations, respectively. The results proved that asking preference information in a parsimonious way is better than requiring an unrealistic cognitive effort to the user. Therefore, this sheds light on the necessity to carefully study how often the user should be queried with a pairwise comparison of solutions to ensure and speed the convergence of the algorithm.

{The paper is structured as follows. In Section \ref{ReviewLocation}, an overview of location problems is provided; MCDA and, in particular, NEMO-II-Ch are presented in Section \ref{MCDA_NEMOIICh}; the particular MOLFP to which we applied NEMO-II-Ch is described in Section \ref{CaseStudy}, while the three algorithms based on the complete knowledge of the user's preference with which NEMO-II-Ch is compared are presented in Section \ref{ComparisonMethods}. The experimental setup and the numerical results are detailed in Section \ref{SetupAndResults}; in Section \ref{DiscussionSection} we discuss the obtained results; finally, the last section provides some conclusions together with possible avenues of research.}

\section{Review on recent approaches to location problems}\label{ReviewLocation}
According to \cite{knuth:art} three types of objectives can be adopted when locating facilities. The $mini-max$ problems, also known as center problems, aim to minimize the maximum distance between a user and its assigned facility \cite{Daskin:1995}. {Several variants of the center problems can be identified (see e.g. \cite{CalikLabbeYaman2015}). For instance, recently, \cite{SchnepperEtAl2019} proposed a new formulation to address a situation  in which the $k$-th largest weighted distance between the users and the facilities needs to be minimized.}  The $mini-sum$ problems minimize the sum of the distances between users and facilities; this objective is well known and much studied and called the median problem \cite{Hakimi:1964,Kariv:1969,Mlad:2007}. Among the median problems let us recall the Discrete Ordered Median Problem (DOMP) \cite{DominguezMarin2013,DominguezMarinEtAl2005,Nickel2001}, where {the objective is the minimization of an ordered weighted average of the distances of the users to the facilities}. Therefore, in this variant of the problem, each user can be seen as an objective. {Lastly, the \textit{covering models} aim to find solutions in which the maximum number of users is covered, i.e. users are positioned within a given threshold distance from a facility \cite{Berman:2010,Church:1974}.}

MOFLPs have been generated from these classical location problems optimizing at the same time more objectives. The very first example was proposed by \cite{Tansel:1982} that optimized the median and the center objectives. Indeed, it proposes to use the median together with other objectives. \cite{Carrizosa:2015} proposed a multiobjective model in which the classical median problem is integrated with a robustness measure that considers potential demand changes. \cite{Blanquero:2002} adopted as additional objective the maximization of the distance from the nearest affected region in order to decrease the impact of the facility on the population. On a similar topic, \cite{Rakas:2004} described a model in which the total number of users that are affected by the facility is minimized. \cite{Pozo:2014} modified the median problem in presence of more DMs considering that the evaluation of the distances between users and facilities is different for every DM. Finally, covering objectives are combined with median objectives as in \cite{Ohsawa:2003}.

Minimizing the distances from the facility (e.g. a disposal site) and the users is also defined in \cite{coutinho2012bi}. They generate solutions containing one of the two objectives adopted (minimizing the distance from the container) and imposing a threshold distance that counts as dissatisfaction from the users.

In addition to the types of objectives described, equality measures can be adopted as objective functions in FLPs \cite{Marsh:1994}. These measures are often used in combination with an efficient objective (e.g., median) to avoid inefficient solutions very much far from all the users \cite{eiseltLaporte:1995}. For example, \cite{Ohsawa:2007} minimized the sum of the absolute differences, the equality measure, and the sum of squared users-facility distances, either to be minimized or maximized for a desirable or obnoxious facility, respectively.

\cite{drezner2006multi} included  all the different types of objectives that have been described so far. They model how to choose the location for a given number of casualty collection points in the California State. They adopt five objectives: the median, the center objective, the covering objectives (using two different distance thresholds) and the variance as equality measure.

It can be noted that many models also include location costs that can depend on several parameters for different potential locations as for example construction costs or maintenance costs \cite{krarup1983simple}. Other MOFLPs adopting several and different objectives can be retrieved in the recent survey by \cite{Farahanisurevey:2010}, often related to the particular case study. They also categorized the MOFLPs on the basis of the methodology developed identifying both exact and heuristic approaches. Beyond that, also several metaheuristics have been applied.
Among these, we focus our attention on the evolutionary algorithms \cite{zhou2011multiobjective}. One first group of applications uses NSGA-II \cite{debieee}. For example,  in \cite{villegas2006solution}  NSGA-II is adopted for the choice of the location of depots in the Colombian coffee supply network maximizing the cover provided by the depots, minimizing the costs of locating the depots and minimizing the distances from purchasing centers to the depot. Similarly, in \cite{bhattacharya2010solving} the NSGA-II methodology is implemented for the location of warehouses and distribution centers in the supply chain perspective, optimizing cost of locating warehouses and cost of transportation from these. Another case for the location of the warehouses in supply chain is reported in \cite{shankar2013location}. In addition to that, some specific applications are approached in \cite{doerner2009multi} for the location of public services in high risk tsunami areas or in \cite{heyns2015multi} for the selection of the best raster points in a Geographical Information System. Finally, \cite{rahmati2014multi} proposed a generic problem in which the first objective function minimizes total setup cost of facilities while the second one minimizes the total expected traveling  and waiting time for the customers.

Other examples of evolutionary algorithms include the application of SPEA2 in \cite{harris2011evolutionary} for deciding the location of depots that serve a single product type to several customers.  Furthermore, the Swarm Optimization has been used as in \cite{yapicioglu2007solving} for approximating the Pareto front in a bi-objective FLP.

While several applications are tackled with evolutionary algorithms, very few examples have been proposed in the literature in which interactive methods have been implemented \cite{Ehrgott:2008}. In \cite{nijkamp1981interactive}, for some generic objective functions of the distances between users and the facility, the DM is asked to indicate some reference levels to be introduced as constraints in the model. Many years later, \cite{karasakal2008interactive} proposed for the two objectives {$mini-max$} and $mini-sum$ an interactive geometrical branch and bound algorithm in which good regions for the location of the facility are selected through the interaction with the DM.
In \cite{dias2008memetic} a memetic algorithm integrates DM's preferences. In particular the DM can choose to indicate reference levels for the objectives, or he can provide the upper bound on the objective function levels. The algorithm can be adapted for several MOFLPs.
A useful tool to help DMs in the interactive phase can be the use of the Geographical Information System (GIS) to help DMs to visualize the potential solutions as in \cite{ALCADAALMEIDA2009111}.
Recently, \cite{DSSCaptivo:2014} developed a Decision Support System for a bi-objective problem; in the computation phase the lexicographic optima and the ideal point  are found, while in the dialog phase  the DM can choose the area in which looking for more non dominated solutions, analyzing maps provided in a GIS environment. This process can be repeated until the DM is satisfied of the final position for the facilities \cite{alves2007review}.

\section{Brief Introduction to MCDA and NEMO-II-Ch}\label{MCDA_NEMOIICh}
\subsection{MCDA and the Choquet integral}\label{MCDAChoquet}
As observed in the previous sections, the use of evolutionary multiobjective optimization methods permits to solve complex multiobjective optimization problems by using evolutionary algorithms. Anyway, the application of these algorithms will give back to the user a set of potentially optimal solutions that will be well-distributed along the Pareto front. The user is therefore asked to choose among them the best one(s) with respect to his preferences. This choice can be very difficult since, in general, the number of non-dominated solutions is quite big and, therefore, the DM could feel himself uncomfortable in performing it.

In order to avoid this, in recent years the interactive methods have been spread out \cite{BDMS_2008}. Their aim is the inclusion of some preference information from the part of the DM addressing the search to the subset of the Pareto front more interesting for him. In order to do that, MCDA methods are used together with evolutionary algorithms (for an updated state of the art survey on MCDA see \cite{GrecoEhrgottFigueira2016}).

Given a set of alternatives $A=\{a,b,\ldots\}$ evaluated on a set of $n$ evaluation criteria $G=\{f_1,\ldots,f_n\}$\footnote{Let us observe that the criteria in MCDA will be the objective functions of the considered multiobjective optimization problem on which the different solutions have to be evaluated.}, MCDA methods deal with ranking, choice and sorting problems. In this case, we will be more interested in ranking and choice problems. In ranking problems, all considered alternatives have to be rank ordered from the best to the worst, while, in choice problems, the best alternative (eventually more than one) has to be chosen, removing all the others. Since the dominance relation\footnote{An alternative $a$ dominates an alternative $b$ iff $a$ is at least as good as $b$ for all considered criteria and better for at least one of them.} stemming from the evaluations of the alternatives on the criteria at hand is too poor, several aggregation methods can be considered. In this paper, we will use as aggregation method the Choquet integral \cite{choquet1953theory} (see \cite{Grabisch1996} for a survey on the use of the Choquet integral in MCDA), a method that can be included under the family of Multiattribute Value Theory (MAVT) \cite{Keeney76}. MAVT methods are based on value functions $U:A\rightarrow\rea$ such that the greater the value assigned to an alternative $a$ by $U$, that is $U(a)$, the better $a$ can be considered. In particular, a preference ($\succ$) and an indifference $(\sim)$ relations can be defined such that $a\succ b$ iff $U(a)>U(b)$, while $a\sim b$ iff $U(a)=U(b)$.

The most common value function $U$ is the additive one

\begin{equation}\label{AdditiveValueFunction}
U(a)=U(f_{1}(a),\ldots,f_n(a))=\sum_{j=1}^{n}u_j(f_j(a)),
\end{equation}

\noindent where, $u_j:A\rightarrow\rea$ are non-decreasing functions of the evaluations $f_j(a)$ for all $f_j\in G.$ Moreover, due to its simplicity, the additive value function most used in applications is the weighted sum

\begin{equation}\label{WeightedSum}
U(a)=U(f_{1}(a),\ldots,f_n(a))=\sum_{j=1}^{n}w_j\cdot f_j(a)
\end{equation}

\noindent where $w_j$ are the weights attached to criteria $f_j\in G$ such that $w_j\geqslant 0$ for all $f_j\in G$ and $\displaystyle\sum_{j=1}^{n}w_j=1$. However, the use of an additive value function assumes that the set of criteria is mutually preferentially independent \cite{Keeney76,wakker1989additive} even if, in real world applications, this assumption is not always verified. Indeed, the evaluation criteria can present a certain degree of positive or negative interaction. On the one hand, two criteria are positively interacting if the importance assigned to them (together) is greater than the sum of the importance assigned to the two criteria taken alone. On the other hand, two criteria are negatively interacting if the importance assigned to them (together) is lower than the sum of the importance assigned to the two criteria singularly. In literature, interaction between criteria is dealt by using non-additive integrals \cite{Grabisch1996,Grabisch2008} and, among them, the most well known is the Choquet integral.

The Choquet integral is based on a capacity, being a set function $\mu:2^{G}\rightarrow\left[0,1\right]$ such that the following constraints are satisfied:
\begin{itemize}
\item[1a)] {$\mu(\emptyset)=0$ and $\mu(G)=1$ (normalization)},
\item[2a)] {$\mu(S)\leqslant\mu(T)$ for all $S\subseteq T\subseteq G$ (monotonicity)}.
\end{itemize}

Given $a\in A$, the Choquet integral of $(f_1(a),\ldots,f_n(a))$ with respect to $\mu$ (in the following, for the sake of simplicity, we shall write ``the Choquet integral of $a$ w.r.t. $\mu$") is computed as follows

\begin{equation}\label{ChoquetCapacity}
C_{\mu}(a)=C_{\mu}(f_{1}(a),\ldots,f_n(a))=\sum_{j=1}^{n}\left[f_{(j)}(a)-f_{(j-1)}(a)\right]\mu(\{f_i\in G: f_i(a)\geqslant f_{(j)}(a)\})
\end{equation}

\noindent where $(\cdot)$ is a permutation of the indices of criteria such that $0=f_{(0)}(a)\leqslant f_{(1)}(a)\leqslant\ldots\leqslant f_{(n)}(a)$\footnote{As observed in \cite{NEMOIICh}, if some evaluations $f_j(a)$ are lower than zero, then it is enough performing a translation $f_{j}(a)\rightarrow f_{j}^{*}(a)=f_j(a)+c$ where $c\geqslant\displaystyle -min_{f_j\in G, a\in A}f_j(a)$ for all $f_j\in G$ and for all $a\in A$.}. To make things easier, a M\"{o}bius transformation of the capacity $\mu$ \cite{Rota,Shafer} and $2$-additive capacities \cite{grabisch1997k} are used in practice:
\begin{itemize}
\item the M\"{o}bius transformation of the capacity $\mu$ is a set function $m:2^G\rightarrow\rea$ such that $\mu(S)=\displaystyle\sum_{T\subseteq S}m(T)$ for all $S\subseteq G$ (conversely, $\displaystyle m(S)=\sum_{T\subseteq S}(-1)^{|S-T|}\mu(T)$ for all $S\subseteq G$) and constraints 1a) and 2a) are replaced by the following ones:
\begin{itemize}
\item[1b)] $m(\emptyset)=0$, $\displaystyle\sum_{T\subseteq G}m(T)=1$,
\item[2b)] for all $f_j\in G$ and for all $S\subseteq G\setminus\{f_j\}$, $\displaystyle\sum_{T\subseteq S}m(T\cup\{f_j\})\geqslant 0$.
\end{itemize}
\noindent In this case, the Choquet integral of $a$ w.r.t. $\mu$ can be written as follows:
\begin{equation}\label{ChoquetMobius}
C_{\mu}(a)=C_{\mu}(f_{1}(a),\ldots,f_n(a))=\displaystyle\sum_{T\subseteq G}m(T)\min_{f_j\in T}f_j(a);
\end{equation}
\item a capacity $\mu$ is said $k$-additive if its M\"{o}bius transformation $m$ is such that $m(T)=0$ for all $T\subseteq G$ such that $|T|>k$. 
\end{itemize}

By using the M\"{o}bius transformation of the capacity $\mu$ and a 2-additive capacity, the Choquet integral can be written in the following linear form

\begin{equation}\label{ChoquetMobius2additive}
C_{\mu}(a)=C_{\mu}(f_{1}(a),\ldots,f_n(a))=\displaystyle\sum_{f_j\in G}m(\{f_j\})f_j(a)+\sum_{\{f_i,f_j\}\subseteq G}m(\{f_i,f_j\})\min\{f_i(a),f_j(a)\}
\end{equation}

\noindent while monotonicity 1b) and normalization constraints 2b) become
\begin{itemize}
\item[1c)] $m(\emptyset)=0$, $\displaystyle\sum_{f_i\in G}m(\{f_i\})+\sum_{\{f_i,f_j\}\subseteq G}m(\{f_i,f_j\})=1$,
\item[2c)] $
\left\{
\begin{array}{l}
m(\{f_i\})\geqslant 0, \;\mbox{for all} \;f_i\in G,\\[1mm]
\displaystyle m(\{f_i\})+\sum_{f_j\in T}m(\{f_i,f_j\})\geqslant 0, \;\mbox{for all} \;f_i\in G \;\mbox{and for all} \;T\subseteq G\setminus\{f_i\}, T\neq\emptyset.\\[1mm]
\end{array}
\right.
$
\end{itemize}

\subsection{NEMO-II-Ch}\label{NEMOIICh}
NEMO-II-Ch \cite{NEMOIICh} is an interactive multiobjective optimization method aiming to address the search to the region of the Pareto front most interesting for the DM. The method belongs to the family of NEMO\footnote{NEMO: Necessary preference enhanced Evolutionary Multiobjective Optimizer} methods \cite{BGSZ-IEEE-TEC-2014} which, on the basis of NSGA-II, integrate some preferences provided by the DM during the iterations of the algorithm. The aim is getting points focused in a particular region of the Pareto front avoiding to waste time in surfing through regions not interesting for the DM. At the beginning, the model uses a simple weighted sum (\ref{WeightedSum}) as preference function and, if necessary, passes to the 2-additive Choquet integral (\ref{ChoquetMobius2additive}) when the preference function is not able to replicate the preferences provided by the DM.

\begin{algorithm*}
\caption{{NEMO-II-Ch} method\label{alg:basic}}
\begin{algorithmic}[1]
\STATE Current preference model = LINEAR.
\STATE Generate initial population of solutions and evaluate them
\REPEAT
\IF{Time to ask the DM}
\STATE Elicit user's preferences by asking DM to compare  two randomly selected non-dominated solutions
\IF{there is no value function remaining compatible with the user's preferences}
\IF{Current preference model = LINEAR}
  \STATE Current preference model = CHOQUET and go to 6:
\ELSE
  \STATE  Remove information on pairwise comparisons, starting from the oldest one, until feasibility is restored and reintroduce them in the reverse order as long as feasibility is maintained
\ENDIF
\ENDIF
\STATE Rank solutions into fronts by iteratively identifying all solutions that are most preferred for at least one compatible value function. Rank within each front using crowding distance
\ENDIF
\STATE Select solutions for mating 
\STATE Generate offspring using crossover and mutation and add them to the population
\STATE Rank solutions into fronts by iteratively identifying all solutions that are most preferred for at least one compatible value function. Rank within each front using crowding distance
\STATE Reduce population size back to initial size by removing worst solutions
\UNTIL{Stopping criterion met}
\end{algorithmic}
\end{algorithm*}

In the following, we shall describe the different steps in Algorithm \ref{alg:basic}:
\begin{itemize}
\item[1:] As mentioned above, at the beginning a linear value function is used to represent the preferences of the DM;
\item[2:] We generate an initial population of solutions and we evaluate them with respect to the considered objective functions;
\item[4-5:] If it is time to ask the DM for preference information, we order the solutions in fronts using the dominance relation, exactly as done in NSGA-II. The non-dominated solutions are put in the first front. Once removed from the population, the other non-dominated solutions are put in the second front and so on, until all solutions have been ordered in different fronts. Inside the same front, the solutions are ordered using the crowding distance \cite{debieee}. The DM is therefore presented with two non-dominated solutions. They are taken in a random way from the first front (if there are at least two solutions in it) or from the following ones having at least two non-dominated solutions. In the extreme case in which there is only one solution for each front (therefore we have a complete order of the solutions), the DM is not presented with any pair of solutions and we can pass to step 15:. \\
Let us suppose that solutions $a$ and $b$ have been chosen to be presented to the DM. He is therefore asked to pairwise compare the two objective functions vectors $(f_1(a),\ldots,f_n(a))$ and $(f_1(b),\ldots,f_n(b))$ stating if $a$ is preferred to $b$ ($a\succ b$), $b$ is preferred to $a$ ($b\succ a$) or $a$ and $b$ are indifferent $(a\sim b)$. A linear constraint will be used to translate this preference information. In particular, $a\succ_{DM} b$ is translated to the constraint $U(a)>U(b)$ and, $a\sim_{DM} b$ iff $U(a)=U(b)$. Let us observe that $U$ is the function in (\ref{WeightedSum}) if the current preference model is the linear one, while $U$ is the function in (\ref{ChoquetMobius2additive}) if the current preference model is the 2-additive Choquet integral;
\item[6:] Checking if there exists at least one value function compatible with the preferences provided by the DM:
\begin{itemize}
\item If the current preference model is the linear one (\ref{WeightedSum}), then one has to solve the following LP problem:

$$
\begin{array}{l}
\;\;\varepsilon^{linear}_{DM}=\max\varepsilon\;\;\mbox{subject to}\\[1mm]
\left.
\begin{array}{l}
U(a)\geqslant U(b)+\varepsilon, \;\mbox{if}\; a\succ_{DM}b,\\[1mm]
{U(a)=U(b), \;\mbox{if}\; a\sim_{DM}b,}\\[1mm]
\displaystyle\sum_{j=1}^{n}w_j=1,\\[1mm]
w_j\geqslant 0,\;\mbox{for all}\;j=1,\ldots,n.
\end{array}
\right\}E^{linear}_{DM}
\end{array}
$$

\noindent Let us observe that one constraint $U(a)\geqslant U(b)+\varepsilon$ should be included for all pairs $(a,b)\in A\times A$ for which the DM states that $a$ is preferred to $b$ ($a\succ_{DM}b$), {while one constraint $U(a)=U(b)$ should be included for all pairs $(a,b)\in A\times A$ for which the DM states that $a$ is indifferent to $b$ ($a\sim_{DM}b$).} If $E^{linear}_{DM}$ is feasible and $\varepsilon^{linear}_{DM}>0$, then there is at least one linear value function compatible with the preferences provided by the DM.
\item If the current preference model is the 2-additive Choquet integral in (\ref{ChoquetMobius2additive}), then one has to solve the following problem:
$$
\begin{array}{l}
\;\;\varepsilon^{Ch}_{DM}=\max\varepsilon\;\mbox{subject to}\\[1mm]
\left.
\begin{array}{l}
C_{\mu}(w_1f_1(a),\ldots,w_nf_n(a))\geqslant C_{\mu}(w_1f_1(b),\ldots,w_nf_n(b))+\varepsilon, \;\mbox{if}\;a\succ_{DM}b,\\[1mm]
C_{\mu}(w_1f_1(a),\ldots,w_nf_n(a))=C_{\mu}(w_1f_1(b),\ldots,w_nf_n(b)), \;\mbox{if}\;a\sim_{DM}b,\\[1mm]
w_j\geqslant 0,\;\mbox{for all}\;j=1,\ldots,n,\\[1mm]
\displaystyle\sum_{j=1}^{n}w_j=1,\\[1mm]
m(\emptyset)=0,\;\mbox{and}\;\displaystyle \sum_{f_i\in G}m(\{f_i\})+\sum_{\{f_i,f_j\}\subseteq G}m(\{f_i,f_j\})=1,\\[1mm]
m(\{f_j\})\geqslant 0, \;\mbox{for all}, \;j=1,\ldots,n,\\[1mm]
m(\{f_{j}\})+\displaystyle\sum_{f_i\in T}m(\{f_i,f_j\})\geqslant 0,\;\mbox{for all}\; j=1,\ldots,n, \;\\
\mbox{and for all}\;T\subseteq\{f_1,\ldots,f_n\}\setminus\{f_j\}, T\neq\emptyset.
\end{array}
\right\}E^{Ch}_{DM}
\end{array}
$$
{Let us underline that in the set of constraints above, we need to introduce a set of weights $\left(w_1,\ldots,w_n\right)$ so that $w_j\geqslant 0$ and $\displaystyle\sum_{j=1}^{n}w_j=1$ since the Choquet integral application implies that all objectives are expressed on the same scale. The set of weights is therefore necessary to put the objectives on the same scale and, for this reason, they become unknown of our model \cite{NEMOIICh}.} \\
If $E^{Ch}_{DM}$ is feasible and $\varepsilon^{Ch}_{DM}>0$, then there is at least one value function, being a 2-additive Choquet integral, compatible with the preferences provided by the DM. Let us observe that the previous problem is not linear anymore and, consequently, we use the Nelder-Mead method \cite{nelder1965simplex} to get the set of weights and the M\"{o}bius parameters optimizing it. It is a numerical algorithm used to solve non-linear optimization problems that, iteratively, {evaluates} solutions belonging to a simplex. At each iteration this simplex is transformed and the procedure continues until a stopping criterion is met (see \cite{NEMOIICh} for a description of the application of the method in this context). The non-linearity of the problem comes from the constraints translating the preferences of the DM since, for all $a\in A$, 
\begin{eqnarray*}
C_{\mu}(w_1f_1(a),\ldots,w_nf_n(a))&=&\displaystyle\sum_{j=1}^{n}w_{j}f_{j}(a)\cdot m\left(\{f_j\}\right) + \sum_{\{f_{i},f_{j}\}\subseteq G}m\left(\{f_{i},f_{j}\}\right)\cdot\min\{w_{i}f_{i}(a),w_{j}f_{j}(a)\}
\end{eqnarray*}
\end{itemize}
and, consequently, $a\succ_{DM}b$ is translated into the constraint 
$$
\displaystyle\sum_{j=1}^{n}w_{j}f_{j}(a)\cdot m\left(\{f_j\}\right)+\sum_{\{f_{i},f_{j}\}\subseteq G}m\left(\{f_{i},f_{j}\}\right)\cdot\min\{w_{i}f_{i}(a),w_{j}f_{j}(a)\}\geqslant 
$$
$$
\displaystyle\sum_{j=1}^{n}w_{j}f_{j}(b)\cdot m\left(\{f_j\}\right)+\sum_{\{f_{i},f_{j}\}\subseteq G}m\left(\{f_{i},f_{j}\}\right)\cdot\min\{w_{i}f_{i}(b),w_{j}f_{j}(b)\}.
$$
Let us underline that in the programming problems above, the strict inequalities have been converted into weak inequalities by using an auxiliary variable $\varepsilon$ which maximization is the objective of our problems. For example, the strict inequality $U(a)>U(b)$ has been converted into the weak inequality $U(a)\geqslant U(b)+\varepsilon$; 
\item[7-10:] If there is not any model compatible with the preferences provided by the DM, we have to distinguish the case in which the current preference model is the linear one from the case in which the current preference model is the 2-additive Choquet integral. In the first case, since {there does not exist} any linear value function able to replicate the preferences of the DM, we increase the complexity of the model passing to the 2-additive Choquet integral. Having more degrees of freedom, it is more flexible and, therefore, it can better adapt itself to the preferences of the DM. In the second case, if we already passed to the 2-additive Choquet integral but there is not any model (therefore weights and M\"{o}bius parameters) compatible with the preferences of the DM, we remove some pieces of this preference information starting from the oldest one until the feasibility is restored. Let us observe that the removal of a piece of preference information should be performed only if the DM agrees on it. This is a relevant aspect since the DM could be very convinced about a certain comparison and, consequently, he doesn't want to remove it;
\item[13:] In order to use the information gathered until now from the DM and, consequently, to address the search to the most interesting region of the Pareto front, we shall order the solutions in fronts in a different way than before. For each solution $x$ in the current population (we shall denote by $A$ the current set of solutions), we have to check if there is at least one compatible function such that $x$ is strictly preferred to all other solutions in $A$. Again, we have to distinguish two cases:
\begin{itemize}
\item If the current preference model is the linear one, the following LP problem has to be solved:
$$
\begin{array}{l}
\;\;\varepsilon_x^{linear}=\max\varepsilon\;\;\mbox{subject to},\\[1mm]
\left.
\begin{array}{l}
U(x)\geqslant U(a)+\varepsilon, \;\mbox{for all}\;a\in A\setminus\{x\},\\[1mm]
E_{DM}^{linear}.
\end{array}
\right\}E_{x}^{linear}
\end{array}
$$
\noindent If $E_{x}^{linear}$ is feasible and $\varepsilon_x^{linear}>0$, then $x$ is put in the first front.
\item If, instead, the current preference model is the 2-additive Choquet integral preference model, then the following programming problem has to be solved:
$$
\begin{array}{l}
\;\;\varepsilon_{x}^{Ch}=\max\varepsilon\;\;\mbox{subject to},\\[1mm]
\left.
\begin{array}{l}
C_{\mu}(w_1f_1(x),\ldots,w_nf_n(x))\geqslant C_{\mu}(w_1f_1(a),\ldots,w_nf_n(a))+\varepsilon, \;\mbox{for all}\;a\in A\setminus\{x\},\\[1mm]
E_{DM}^{Ch}.
\end{array}
\right\}E_{x}^{Ch}
\end{array}
$$
If $E_{x}^{Ch}$ is feasible and $\varepsilon_{x}^{Ch}>0$, then $x$ is put in the first front.
\end{itemize}
Once the first front has been built, all solutions contained in it are removed from the current population and the same procedure is used with the remaining solutions to build the second front. We shall continue in this way until all solutions have been ordered in different fronts. Inside the same front, solutions are ordered using the crowding distance. \\
{In the rare case in which there is not any solution that can be preferred to the others for any compatible model, all solutions are retained equally preferable and, therefore, they are put in the same front;}
\item[15-18:] The usual evolution of the population is performed by using the selection, crossover and mutation operators together with the ordering of the population described above;
\item[3-19:] Repeat steps 4-18 until the stopping condition has been met.
\end{itemize}

\section{Using Interactive Evolutionary Multiobjective Optimization in location problems: a case study}\label{CaseStudy}
We test our approach on a well known multiobjective location problem introduced in \cite{drezner2004location} and later in \cite{drezner2006multi}. The problem, considered {as} a reference in its domain, consists in choosing the location of a given number $p$ of facilities among a set of potential locations optimizing five different classical objective functions for FLPs. More in detail, the facilities are Casualty Collection Points (CCPs) to which people can go if they need help in case disasters have happened. These centers should operate where a huge amount of people need to be provided with emergency service. In \cite{drezner2004location} a comparison of the different objectives is proposed and also a first multiobjective version, including only three objectives, is formulated; whereas in \cite{drezner2006multi} a multiobjective heuristic has been introduced adopting the five objective functions described later.
The problem is of particular interest among the MOFLPs because  at least one $mini-max$ objective is selected, one for the $mini-sum$, and one equality measure are optimized.

We define:
\begin{itemize}
  \item $I=\{1,\ldots , q\}$: the set of demand points,
  \item $L=\{1,\ldots , m\}$: the set of potential locations for the facilities,
  \item $d_{ij}$: the distance between demand point $i$ and potential facility $j$,
  \item $pop_i$: the population at the demand point $i$,
  \item $p$: the total number of facilities to locate,
  \item $P\subseteq L$: a vector of $p$ selected facilities in $L$,
  \item $D_i(P)$: the distance from a demand point $i$ to the closest facility in $P$,
      $$D_i(P)=\min_{k\in P} \{d_{ik}\}.$$
\end{itemize}

We consider five objectives:

\begin{enumerate}
\item  \textit{The median objective}, minimizes the sum of the distances between the demand points and the closest facility \cite{Hakimi:1964,revelleswaim:1970}:
$$\min_{P}f_1(P)=\min_{P}\left[\frac{1}{q}\sum_{i=1}^{q}D_i(P)\right],$$
\item  \textit{The maximum distance objective}, minimizes the distance of the farthest demand point \cite{Hakimi:1964}:
    $$\min_{P}f_2(P)=\min_P\{\max_i\{D_i(P)\}\},$$
\item  \textit{The maximum covering objectives}, maximize {the population} inside  two different distance thresholds  $S_1$ and $S_2$ \cite{Church:1974}:
$$\max_{P}f_3(P)=\max_{P}\sum_{i:\;D_i(P)\leqslant S_1}pop_i,$$
$$\max_{P}f_4(P)=\max_P\sum_{i:\;D_i(P)\leqslant S_2}pop_i.$$
\item  \textit{The minimum variance objective}, balances the distances between demand points and the closest facility, minimizing the variance of the closest distances for all the demand points \cite{Marsh:1994}:
$$\min_P f_5(P)=\min_P\displaystyle\frac{\sum_{i=1}^{q}[D_i(P)-f_1(P)]^2}{q}.$$
\end{enumerate}

The case study has ${I=\{1,\ldots , 577\}}$ demand points and ${L=\{1,\ldots , 141\}}$ potential sites for the facilities located in the Orange County in California, an area where the careful planning for the location of CCPs  represents an essential requirement due to  frequent earthquakes. The data, that include  coordinates  and associated weights for the demand points and coordinates for the potential facilities, are available upon request to the authors of \cite{drezner2006multi}.

Let us point out that this is just one of the possible examples that our methodology can handle. Our approach is very flexible and we could adopt many different objective functions.

{
\section{Algorithms used for the comparison}\label{ComparisonMethods}
As already observed above, the use of a heuristic not taking into account preferences of the DM, such as NSGA-II, gives back the user with a set of non-dominated vectors of $p$-facilities from which he has to choose the best with respect to his preferences. For this reason, we proposed to apply NEMO-II-Ch to address the search not to the entire Pareto front but to the most interesting part for the user. \\
We shall consider the full size problem in which the 141 different locations will be taken into account choosing the best $p$ among them with $p=4,5.$ Moreover, we will simulate different users' value functions. On the one hand, we will show that, in most of the cases, NEMO-II-Ch is able to find the best subset of $p$ locations for the user by asking few preference information. On the other hand, to test its performances, we will compare them to the performances of three algorithms, denoted by EA-UVF, EA-UVF1 and EA-UVF2. These are based on the knowledge of the user's true value function that is, instead, unknown to the NEMO-II-Ch algorithm. While the EA-UVF algorithm has been presented in \cite{NEMOIICh}, its two variants, namely EA-UVF1 and EA-UVF2, are presented in this paper for the first time. The three algorithms are briefly presented in the following sections. }

{
\subsection{\textit{EA-UVF: Evolutionary Algorithm based on User's Value Function}}
This algorithm has been presented in \cite{NEMOIICh} and its main steps, which are listed in Algorithm \ref{EAUVFAlgorithm} are detailed in the following lines: 
}

{
\begin{algorithm*}
\caption{Evolutionary Algorithm User's Value Function (EA-UVF) algorithm\label{EAUVFAlgorithm}}
\begin{algorithmic}[1]
\STATE Generate initial population of solutions and evaluate them
\STATE Compute the utility of each solution by the user's true value function
\STATE Rank the solutions into fronts with respect to their true value
\REPEAT
\STATE Select solutions for mating
\STATE Generate offspring using crossover and mutation and add them to the population
\STATE Rank the solutions into fronts with respect to their true value
\STATE Reduce population size back to initial size by removing worst solutions
\UNTIL{Stopping criterion met}
\end{algorithmic}
\end{algorithm*}
}

{
\begin{itemize}
\item[1:] Generate an initial population of solutions and evaluate them with respect to the considered objective functions;
\item[2:] Compute the utility of each solution by using the user's true value function;
\item[3:] Rank the solutions into fronts by using the values assigned to them from the user's true value function and computed at the previous step. The solution having the best utility value (the minimum [maximum] value if the user's true value function has to be minimized [maximized]) is put in the first front; the solution having the second best utility value is put in the second front and so on until the solution having the worst utility value that is included in the last front. Solutions having the same utility value are included in the same front;
\item[5:8] Evolve the population; 
\item[4-9:] Repeat steps 5-8 until the stopping condition has not been met.\\
\end{itemize}
}

{
\subsection{\textit{EA-UVF1: NSGA-II with diversification replaced by User's Value Function}}
The steps of the EA-UVF1 algorithm are shown in Algorithm \ref{EAUVFAlgorithm1} and detailed in the following lines:
}

{
\begin{algorithm*}
\caption{NSGA-II with diversification replaced by User's Value Function (EA-UVF1) \label{EAUVFAlgorithm1}}
\begin{algorithmic}[1]
\STATE Generate initial population of solutions and evaluate them
\STATE Rank solutions into fronts by dominance and inside each front order them using their true value
\REPEAT
\STATE Select solutions for mating
\STATE Generate offspring using crossover and mutation and add them to the population
\STATE Rank solutions into fronts by dominance and inside each front order them using their true value
\STATE Reduce population size back to initial size by removing worst solutions
\UNTIL{Stopping criterion met}
\end{algorithmic}
\end{algorithm*}
}

{
\begin{itemize}
\item[1:] Generate an initial population of solutions and evaluate them with respect to the considered objective functions;
\item[2:] Rank solutions in non-dominated fronts. Then, inside each front, compute the true value of all solutions and rank them by these utility values;
\item[4-7:] Evolve the population;
\item[3-8:] Repeat steps 4-7 until the stopping condition has not been met.
\end{itemize}
}

{
\noindent The EA-UVF1 implements exactly the NSGA-II method with the replacement of the crowding distance used to diversify solutions inside the same front with the value assigned to the solutions by the user's true value function. 
}

{
\subsection{\textit{EA-UVF2: NSGA-II with a roulette wheel driven by User's Value Function}}
The steps of the EA-UVF2 algorithm are shown in Algorithm \ref{EAUVFAlgorithm2} and detailed in the following lines:
}

{
\begin{algorithm*}
\caption{NSGA-II with a roulette wheel driven by User's Value Function (EA-UVF2) \label{EAUVFAlgorithm2}}
\begin{algorithmic}[1]
\STATE Generate initial population of solutions and evaluate them
\REPEAT
\STATE Assign a probability to be parent to each solutions by using their true value
\STATE Select solutions for mating
\STATE Generate offspring using crossover and mutation and add them to the population
\STATE Rank solutions into fronts by dominance and inside each front order them by the crowding distance
\STATE Reduce population size back to initial size by removing worst solutions
\UNTIL{Stopping criterion met}
\end{algorithmic}
\end{algorithm*}
}

{
\begin{itemize}
\item[1:] Generate an initial population of solutions and evaluate them with respect to the considered objective functions;
\item[3:] A probability to be parent of the next generation is assigned to each solution in the population. This probability, denoted by $Prob(P)$, is computed as 
\begin{equation}\label{Prob1}
Prob(P)=\frac{U(P)}{\displaystyle\sum_{P\in POP}U(P)} \qquad\mbox{if $U$ has to be maximized},
\end{equation}
\begin{equation}\label{Prob2}
Prob(P)=\frac{\frac{1}{U(P)}}{\displaystyle\sum_{P\in POP}\frac{1}{U(P)}} \qquad\mbox{if $U$ has to be minimized}
\end{equation}
\noindent and $POP$ denotes the current population of solutions;
\item[4-7:] Evolve the population;
\item[2-8:] Repeat steps 3-7 until the stopping condition has not been met. \\
\end{itemize}
The EA-UVF2 algorithm implements, therefore, all steps of the NSGA-II method, while the user's true value function is used to assign a probability to be parent of the next generation to each solution. The better the value assigned by the user's true value function to a solution, the higher its probability to become parent of the next generation. \\
}

{
Let us conclude this section by underlining that the EA-UVF represents the ideal situation in which the algorithm knows exactly how the user chooses among two whichever solutions and, consequently, it has the maximal theoretical availability of preference information. At the same time, the EA-UVF1 and the EA-UVF2 use this information, on the one hand, to select solutions within non-dominated fronts the generated population and, on the other hand, to decide which solutions are the best to be parents of the next generation. However, all of them use the whole preference information that the DM could theoretically provide by preferentially ranking all solutions at all iterations of the evolutionary algorithm. Of course, in real life the DM can not be able to provide all these preferences because of the unrealistic huge cognitive burden related to the request of so many preference comparisons at each iteration. For this reason, a methodology being much more parsimonious in asking preferences to the DM is requested for any real world application. To study the amount of preference information necessary to get reasonably acceptable solutions, we investigate the relation between, on the one hand, the frequency of asking preferences to the user and, on the other hand, the quality of results and the speed of the algorithms' convergence. To this aim, in the following simulations, we run NEMO-II-Ch asking the DM one preference every 5, 10 and 20 generations, respectively.
}

\section{Experimental setup and numerical results}\label{SetupAndResults}
The parameters and the technical details used in the simulations are the following:
\begin{itemize}
\item The population $POP$ is composed of 30 solutions where each solution is a vector $P$ of $p$ different integer values taken in the interval $[1,m]$;
\item The mating selection is performed by tournament selection in all methods apart from EA-UVF2 where it is performed by a roulette wheel selection:
\begin{itemize}
\item 
{
\textit{Tournament selection:} Let us denote by $P_{1},\ldots,P_{30}$ the solutions in the current population. To each solution $P_{s}$ is associated the front it belongs to ($F_s$). Moreover, in all methods each solution is associated with a \textit{second score}. In NEMO-II-Ch and in EA-UVF2 this second score is the crowding distance ($CD_s$)\footnote{Citing \cite{debieee}, the crowding distance is \textit{...``the average distance of two points on either side of a particular solution along each of the objectives}" and it is computed to maintain the diversification of the population. The higher the crowding distance of a solution $P_{s}$, the more isolated is the solution in the considered population.}, while, in EA-UVF1 the second score is the true value. We create a random permutation of the solutions in the population denoted by $P_{(1)},\ldots,P_{(30)}$. Then, a tournament is performed between $P_{s}$ and $P_{(s)}$ for each $s=1,\ldots,30,$ to choose which solution has to be selected as parent of the next generation. The tournament is won from the solution being in the lowest front ($P_{s}$ iff $F_{s}<F_{(s)}$ or $P_{(s)}$ iff $F_{(s)}<F_{s}$) or, if they belong to the same front ($F_{s}=F_{(s)}$), from the solution having the greatest second score. If $P_{s}$ and $P_{(s)}$ belong to the same front and they have the same second score, the winner is chosen randomly. 30 tournaments will therefore be performed and, consequently, 30 solutions will become parents of the next generation. Denoting by $P_{s}^{'}$ the winner of the tournament between $P_{s}$ and $P_{(s)}$, the pairs of parents which will generate the offsprings of the next generation are, therefore, $(P^{'}_{1},P^{'}_{2})$, $(P^{'}_{3},P^{'}_{4})$,$\ldots$,$(P^{'}_{29},P^{'}_{30})$;
\item \textit{Roulette wheel selection:} Since, as in the tournament selection, 15 pairs of parents $(P^{'}_{1},P^{'}_{2})$, $(P^{'}_{3},P^{'}_{4})$,$\ldots$,$(P^{'}_{29},P^{'}_{30})$ have to be chosen, for each $k=1,\ldots,30$, a solution is sampled in a random way from the probability distribution given by eq. (\ref{Prob1}) if the user's true value function $U$ has to be maximized or by eq. (\ref{Prob2}) if the same function as instead to be minimized; the sampled solution becomes, therefore, the parent $P_{k}^{'}$ of the next generation; 
}
\end{itemize} 
\item Each pair of parents generate two offspring by one-point crossover with probability of 1 and random resetting mutation\footnote{\textit{``...in each position independently, with probability $p_m$, a new value is chosen at random from the set of permissible values"} \cite{Eibe03a}} with probability of $\frac{1}{p}$ \cite{Eibe03a}; in particular, since each solution can contain a certain location at most once, the one-point crossover has to be slightly modified if the two considered solutions have some common locations. In this case, the common potential location(s) are inherited by both offspring, while the one-point crossover is performed on the two vectors composed by uncommon potential locations for both parents. For example, let us suppose that the two parents solutions are (10,15,21,30) and (6,10,20,50). In this case, the potential location labeled by 10 is present in both parents and, therefore, it is inherited by the two offspring. The remaining vectors composed of uncommon locations are (15,21,30) and (6,20,50). To these two vectors the one-point crossover is applied exchanging the two tails. Supposing that the cut point is the second integer, exchanging the two tails we obtain the vectors (15,21,50) and (6,20,30). The two offspring will therefore be the vectors (10,15,21,50) and (6,10,20,30). \\
Let us underline that the evolution of the population is performed in such a way that if a new offspring is exactly the same as another solution in the current population, it is ``killed". Therefore, it is not possible having multiple copies of the same solutions in the population;
\item Considering the set $L$ of potential locations and a solution $P$ composed of $p$ of these potential locations, we assumed the following different forms of user's preferences described as follows:
\begin{itemize}
\item[$U^D$)] the maximal deviation from the optimal objective values \cite{drezner2006multi} is computed as follows 
$$
U^{D}(P)=\displaystyle\max_{k\in\{1,\ldots,5\}}\left\{\Delta_k(P)\right\}
$$ 
\noindent where
$$
\Delta_{k}(P)=\left\{
\begin{array}{ll}
\frac{f_k(P)-f_{k}^{*}}{f_{k}^{*}}, & \mbox{if the objective $f_k$ is to be minimized},\\[3mm]
\frac{f_{k}^{*}-f_k(P)}{f_{k}^{*}}, & \mbox{if the objective $f_k$ is to be maximized},\\
\end{array}
\right.
$$
\noindent and
$$
f_{k}^{*}=\left\{
\begin{array}{ll}
\displaystyle f_{k}^{min}=\min_{\overline{P}\subseteq L: \;|\overline{P}|=p}f_{k}(\overline{P}), & \;\mbox{if the objective $f_k$ is to be minimized},\\
\displaystyle f_{k}^{max}=\max_{\overline{P}\subseteq L: \;|\overline{P}|=p}f_{k}(\overline{P}), & \;\mbox{if the objective $f_k$ is to be maximized},\\
\end{array}
\right.
$$
\noindent that is, $f_k^{*}$ is the optimal value for the objective $f_k$, $k=1,\ldots,5;$ a solution $P$ is preferred to a solution $P'$ if $U^D(P) < U^D(P')$;
\item[$U^{D}_{v})$]  On the basis of the $U^{D}$ {defined} above, we considered the function $U^{D}_{v}$ computed as follows: 
$$
U^{D}_{v}(P)=\max_{k\in v}\left\{\Delta_{k}(P)\right\}
$$
\noindent where $v\in\left\{\{1,2,3,4\},\{1,2,3,5\},\{1,2,4,5\},\{1,3,4,5\},\{2,3,4,5\}\right\}$. In this way, we shall take into account only four of the five objective functions simultaneously;
\item[$U^N)$] the value is computed {as follows}
$$
U^{N}(P)=\displaystyle\sum_{k=1}^{5}w_k\cdot\overline{f}_{k}(P)
$$
\noindent where 
$$
\overline{f}_{k}(P)=
\left\{
\begin{array}{ll}
\frac{f_{k}(P)-f_{k}^{min}}{f_{k}^{max}-f_{k}^{min}}, & \mbox{if the objective $f_k$ is to be minimized},\\[0,25cm]
\frac{f_{k}^{max}-f_{k}(P)}{f_{k}^{max}-f_{k}^{min}}, & \mbox{if the objective $f_k$ is to be maximized},
\end{array}
\right.
$$
\noindent $w=(0.1,0.15,0.2,0.25,0.3)$, and a solution $P$ is preferred to a solution $P'$ if $U^N(P) < U^N(P')$;
\item[$U_{v}^N)$] the value is computed as follows
$$
U_{v}^N(P)=\displaystyle\sum_{k\in v}w'_{k}\cdot\overline{f}_{k}(P)
$$
\noindent where $w'=\left(0.1,0.2,0.3,0.4\right)$ and $v\in\left\{\{1,2,3,4\},\{1,2,3,5\},\{1,2,4,5\},\{1,3,4,5\},\{2,3,4,5\}\right\}$\footnote{Let us observe that in the computation of $U_{v}^{N}(P)$ the functions $\overline{f}_k$ have a weight increasing with $k$. For example, if $v=\{1,3,4,5\}$, then, $U_{v}^{N}(P)=0.1\cdot \overline{f}_{1}(P)+0.2\cdot \overline{f}_{3}(P)+0.3\cdot \overline{f}_{4}(P)+0.4\cdot \overline{f}_{5}(P)$.}. Also in this case we consider a subset composed of four of the five objective functions and a solution $P$ is preferred to a solution $P'$ if $U_{v}^N(P) < U_{v}^N(P')$.
\end{itemize}
\noindent For all considered user's value functions, the best subset of $p$ locations is $P_b\subseteq L$, such that $|P_b|=p$ and $U(P_{b})=\displaystyle\min_{\overline{P}\subseteq L:\; |\overline{P}|=p}U(\overline{P})$ where $U\in\{U^D,U^{D}_{v},U^{N},U^{N}_{v}\}$; 
\item All algorithms are run for a maximum of 1,000 generations. In particular, for NEMO-II-Ch we asked the user to provide one preference comparison every 5, 10 and 20 generations. The resulting algorithms are therefore denoted by NIICh\_5, NIICh\_10 and NIICh\_20. All the algorithms stop as soon as $P_{b}$ is present in the current population or when the maximum number of generations has been reached.
\end{itemize}

After we described the setup of the simulations, let us present the results of the application of the compared methods to the considered full-size problem. This means that we shall check for the best subset of $p$ locations, with $p=4,5$, among the 141 taken into account. Of course, this problem is quite difficult since the possible subsets of $p$ locations from which the best has to be discovered are $\binom{141}{4}=15,777,195$ and $\binom{141}{5}=432,295,143$, respectively. Therefore, we would like to prove that the method is able to deal with big-size problems in which a huge number of solutions is involved. We performed 50 independent runs (changing, therefore, the starting population), and we applied the three NEMO-II-Ch variants (NIICh\_5, NIICh\_10 and NIICh\_20) as well as the three algorithms knowing the user's true value function (EA-UVF, EA-UVF1 and EA-UVF2).\\ 
In the tables below, to present the results of the performed simulations, we used the following notation: 
\begin{itemize}
\item $\#SR$: number of runs (over the 50 considered), in which the algorithm was able to discover the best subset $P_b$ of possible locations; 
\item $M\#G$: mean number of generations necessary to the algorithm to discover $P_{b}$;
\item $S\#G$: standard deviation of the number of generations necessary to the algorithm to discover $P_b$;
\item {$A\#P$: mean number of pairwise comparisons asked to the user necessary to discover $P_{b}$. We did not include this data for EA-UVF and EA-UVF1 since they are only used as benchmark and a comparison between the number of pairwise comparisons asked from the NEMO-II-Ch versions and the one involved in the application of both algorithms is meaningless. Of course, the number of times the user is queried by NEMO-II-Ch is only a small portion of the number of times the user has to provide a pairwise comparison in the two algorithms. Just to give an example, let us underline that in EA-UVF and EA-UVF1, where solutions are ranked with respect to the user's true value function, to rank order $p$ solutions it is necessary to perform $\frac{p(p-1)}{2}$ pairwise comparisons\footnote{The best solution is found after $p-1$ comparisons, the second after $p-2$ comparisons and so on.}. This means that to rank order 30 solutions in the population, the user has to provide 435 pairwise comparisons in a single iteration and, as will be clear in the next section, this number is much higher than the number of pairwise comparisons asked from the three NEMO-II-Ch versions in whichever considered test problem. \\
With respect to EA-UVF2, the user is not asked to provide any pairwise comparison. However, the algorithm can never be applied in practice since it is assumed that the user is able to assign a utility to each solution, utility that needs to be used to implement the roulette wheel selection described above. Of course, this is not realistic at all;}
\item $S\#P$: standard deviation of the number of pairwise comparisons asked to the user necessary to discover $P_{b}$;
\item $MT$: mean time necessary to the algorithm to discover $P_{b}$; all simulations have been performed using the commercial software MATLAB2019 but on different PCs. For each method and each user's value function, the 50 runs have been performed on the same machine. In the tables presenting the results, we reported the characteristics of the PCs used to perform the different simulations; 
\item $ST$: standard deviation of the time necessary to the algorithm to discover $P_{b}$\footnote{The mean and the standard deviation are computed for the runs in which $P_b$ is discovered.};  
\item $A\_BRSD$: this is the average distance of the best solution in the final population from the optimal solution $P_{b}$.The distance, denoted by $BRSD(U)$, is computed only for the simulations in which the algorithm was not able to discover $P_{b}$ (in the case in which the algorithm is able to discover $P_{b}$ the distance is zero). Denoting by $P^{Best}$ the best solution in the final population, following \cite{TomczykKadzinski2019}, $BRSD(U)$ is computed as

\begin{equation}\label{BRSD}
BRSD(U)=\frac{|U(P^{Best})-U(P_b)|}{U(P_{b})}. 
\end{equation}
 
The less $BRSD(U)$, the better the performance of the algorithm. The value $A\_BRSD$ is then obtained by averaging $BRSD(U)$ over the number of runs in which the algorithm was not able to discover $P_{b}$.
\end{itemize}

\subsection{Comparison with EA-UVF, EA-UVF1 and EA-UVF2}

\renewcommand\arraystretch{1.3}

{
\begin{table}[!h]
	\centering
\caption{Results for functions $U^{N}$ and $U_{v}^N$ considering $p=4$. All simulations have been performed with four {different} PCs which characteristics and labels are the following: (PC1) intel core i7 3.6GHz; (PC2) intel core i5 2.5GHz; (PC3) intel core i7 2.7GHz; (PC4) intel core i7 1.9GHz. \label{N4Table}}
\resizebox{0.7\textwidth}{!}{
		\begin{tabular}{ccccccc}
		\hline
    $U^{N}$ & \textbf{NIICh\_5} (PC3) & \textbf{NIICh\_10} (PC3) & \textbf{NIICh\_20} (PC3) & \textbf{EA-UVF} & \textbf{EA-UVF1} & \textbf{EA-UVF2} \\
		\hline
		$\#SR$ & \textbf{50/50}  & \textbf{50/50}  & \textbf{50/50}  & \textbf{50/50}  & \textbf{50/50} & 49/50  \\
		$M\#G$ & 80.92  & 80.44  & 113.16 & 79.7   & 77.54 & 152.78 \\ 
		$S\#G$ & 55.28  & 46.99  & 76.45  & 50.34  & 49.16 & 137.43 \\
		$A\#P$ & 16.80  & 8.60   & \textit{6.16}   &   & & \\
		$S\#P$ & 11.04  & 4.69   & 3.83   &   & & \\ 
		$MT$   & 51.71s & 36.4s & 46.87s  &  & &  \\ 
		$ST$   & 43.99s & 24.29s & 33.16s  &  & &  \\  
		$A\_BRSD$ &        &  &    &  & &  0.406   \\
		\hline
		$U_{1234}^{N}$ & \textbf{NIICh\_5} (PC3) & \textbf{NIICh\_10} (PC2) & \textbf{NIICh\_20} (PC2) & \textbf{EA-UVF} & \textbf{EA-UVF1} & \textbf{EA-UVF2} \\
		\hline
		$\#SR$ & \textbf{50/50}  & \textbf{50/50}  & \textbf{50/50}  & \textbf{50/50} & \textbf{50/50} & \textbf{50/50}  \\
		$M\#G$ & 73.98  & 74.04  & 85.06  & 65.6  & 65.90 & 134.54  \\ 
		$S\#G$ & 61.49  & 40.43  & 72.10  & 51.15 & 49.1  & 153.43  \\
		$A\#P$ & 15.42  & 7.98   & \textit{4.78}   &  & &  \\
		$S\#P$ & 12.27  & 4.00   & 3.65   &  & & \\ 
		$MT$   & 57.05s & 1.85m  & 56.79s &  & &   \\ 
		$ST$   & 1.13m  & 1.41m  & 52.01s &  & &  \\
		\hline
		$U_{1235}^{N}$ & \textbf{NIICh\_5} (PC2) & \textbf{NIICh\_10} (PC2) & \textbf{NIICh\_20} (PC2) & \textbf{EA-UVF} & \textbf{EA-UVF1} & \textbf{EA-UVF2} \\
		\hline
		$\#SR$ & \textbf{50/50}  & \textbf{50/50}  & \textbf{50/50}  & \textbf{50/50} & \textbf{50/50} & \textbf{50/50}  \\
		$M\#G$ & 140.88 & 138.02 & 152.5  & 86.66 & 116.64 & 225.54  \\ 
		$S\#G$ & 100.30 & 106.56 & 109.93 & 68.12 & 95.84 & 196.90   \\
		$A\#P$ & 28.76  & 14.34  & \textit{8.06}   &  & & \\
		$S\#P$ & 20.05  & 10.64  & 5.58   &  & & \\ 
		$MT$   & 3.70m  & 2.83m  & 2.16m  &  & &  \\ 
		$ST$   & 4.04m  & 2.98m  & 2.05m  &  & &  \\
		\hline
		$U_{1245}^{N}$ & \textbf{NIICh\_5} (PC3) & \textbf{NIICh\_10} (PC4) & \textbf{NIICh\_20} (PC1) & \textbf{EA-UVF} & \textbf{EA-UVF1} & \textbf{EA-UVF2} \\
		\hline
		$\#SR$ & \textbf{50/50}  & \textbf{50/50}  & \textbf{50/50}  & \textbf{50/50} & \textbf{50/50} & 5/50  \\
		$M\#G$ & 122.66 & 143.56 & 189.36 & 89.86 & 117.34 & 117.2  \\ 
		$S\#G$ & 61.13  & 88.49  & 103.44 & 77.42 & 73.69 & 63.14  \\
		$A\#P$ & 25.12  & 14.92  & \textit{10.02}  &  & & \\
		$S\#P$ & 12.24  & 8.81   & 5.20   &  & & \\ 
		$MT$   & 2.1m   & 1.9m   & 2.16m  &  & &   \\ 
		$ST$   & 1.47m  & 1.9m   & 1.5m   &  & &  \\
		$A\_BRSD$ &        &  &    &  & & 0.233     \\
		\hline
		$U_{1345}^{N}$ & \textbf{NIICh\_5} (PC4) & \textbf{NIICh\_10} (PC4) & \textbf{NIICh\_20} (PC4) & \textbf{EA-UVF} & \textbf{EA-UVF1} & \textbf{EA-UVF2} \\
		\hline
		$\#SR$ & \textbf{50/50}  & \textbf{50/50}  & \textbf{50/50}  & \textbf{50/50} & \textbf{50/50} & \textbf{50/50}  \\
		$M\#G$ & 76.92  & 82.86  & 112.22 & 63.08 & 76.48 & 162.52  \\ 
		$S\#G$ & 51.36  & 48.87  & 79.03  & 40.21 & 60.42 & 150.01  \\
		$A\#P$ & 16.08  & 8.88   & \textit{6.14}   &  & &  \\
		$S\#P$ & 10.26  & 4.90   & 3.99   &  & &  \\ 
		$MT$   & 1.25m  & 59.62s & 1.14m  &  & &  \\ 
		$ST$   & 1.07m  & 43.73s & 52.14s &  & &  \\
		\hline
		$U_{2345}^{N}$ & \textbf{NIICh\_5} (PC4)& \textbf{NIICh\_10} (PC4) & \textbf{NIICh\_20} (PC4) & \textbf{EA-UVF} & \textbf{EA-UVF1} & \textbf{EA-UVF2} \\
		\hline
		$\#SR$ & \textbf{50/50}  & \textbf{50/50}  & \textbf{50/50}  & \textbf{50/50} & \textbf{50/50} & \textbf{50/50}  \\
		$M\#G$ & 77.46  & 83.02  & 108.2  & 68.86 & 84.56 & 119.76  \\ 
		$S\#G$ & 53.44  & 51.36  & 72.09  & 45.84 & 62.00 & 97.92  \\
		$A\#P$ & 16.16  & 8.86   & \textit{5.94}   &  & &  \\
		$S\#P$ & 10.67  & 5.12   & 3.61   &  & &  \\ 
		$MT$   & 1.26m  & 59.47s & 1.04m  &  & &  \\ 
		$ST$   & 1.18m  & 48.9s  & 43.31s &  & &  \\
		\end{tabular}		
		}		
\end{table}
}

{
\begin{table}[!h]
	\centering
\caption{Results for functions $U^{N}$ and $U_{v}^N$ considering $p=5$. All simulations have been performed with four {different} PCs which characteristics and labels are the following: (PC1) intel core i7 3.6GHz; (PC2) intel core i5 2.5GHz; (PC3) intel core i7 2.7GHz; (PC4) intel core i7 1.9GHz.\label{N5Table}}
\resizebox{0.7\textwidth}{!}{
		\begin{tabular}{ccccccc}
		\hline
    $U^{N}$ & \textbf{NIICh\_5} (PC2) & \textbf{NIICh\_10} (PC2) & \textbf{NIICh\_20} (PC2) & \textbf{EA-UVF} & \textbf{EA-UVF1} & \textbf{EA-UVF2} \\
		\hline
		$\#SR$ & \textbf{50/50}  & \textbf{50/50}  & \textbf{50/50}  & \textbf{50/50}  & \textbf{50/50} & 49/50   \\
		$M\#G$ & 140.94 & 156.08 & 184.96 & 112.22 & 137.20 & 269.08 \\ 
		$S\#G$ & 68.73  & 66.12  & 94.43  & 69.82  & 86.82 & 192.38 \\
		$A\#P$ & 28.82  & 16.16  & \textit{9.78}  &  & & \\
		$S\#P$ & 13.75  & 6.63   & 4.68   &  & & \\ 
		$MT$   & 2.25m  & 1.92m  & 1.98m  &  & &  \\ 
		$ST$   & 1.35m  & 1.05m  & 1.19m  &  & &  \\
		$A\_BRSD$ &        &  &    &  & & 0.127    \\
		\hline
		$U_{1234}^{N}$ & \textbf{NIICh\_5} (PC4) & \textbf{NIICh\_10} (PC4) & \textbf{NIICh\_20} (PC4) & \textbf{EA-UVF}  & \textbf{EA-UVF1} & \textbf{EA-UVF2} \\
		\hline
		$\#SR$ & \textbf{50/50}  & \textbf{50/50}  & \textbf{50/50}  & \textbf{50/50} & \textbf{50/50} & 45/50  \\ 
		$M\#G$ & 110.64 & 143.32 & 166.54 & 93.96 & 110.74 & 297.62  \\ 
		$S\#G$ & 56.28  & 64.70  & 88.90  & 59.80 & 52.4 & 230.45  \\
		$A\#P$ & 22.74  & 14.90  & \textit{8.86}   &  & &  \\
		$S\#P$ & 11.22  & 6.45   & 4.43   &  & &  \\ 
		$MT$   & 2.08m  & 2.59m  & 2.01m  &  & &  \\ 
		$ST$   & 1.36m  & 1.51m  & 1.48m  &  & &  \\		
		$A\_BRSD$ &        &  &    &  & & 0.283     \\
		\hline
		$U_{1235}^{N}$ & \textbf{NIICh\_5} (PC4) & \textbf{NIICh\_10} (PC4) & \textbf{NIICh\_20} (PC4) & \textbf{EA-UVF} & \textbf{EA-UVF1} & \textbf{EA-UVF2} \\
		\hline
		$\#SR$ & \textbf{50/50}  & 49/50  & 47/50  & 45/50   & \textbf{50/50}   & 11/50    \\ 
		$M\#G$ & 313.8  & 308.31 & 313.96 & 354.22  & 256.6  & 638.27 \\ 
		$S\#G$ & 225.6  & 214.84 & 209.36 & 228.32  & 167.76 & 238.52  \\
		$A\#P$ & 63.44  & 31.35  & \textit{16.30}  &   & & \\
		$S\#P$ & 45.14  & 21.44  & 10.48  &   & & \\ 
		$MT$   & 11.38m & 9.13m  & 9.04m  &   & &  \\ 
		$ST$   & 12.13m & 9.28m  & 7.63m  &   & &   \\		
		$A\_BRSD$ &     & 0.007  & 0.007  & 0.007 & & 0.213    \\
		\hline
		$U_{1245}^{N}$ & \textbf{NIICh\_5} (PC4) & \textbf{NIICh\_10} (PC4) & \textbf{NIICh\_20} (PC4) & \textbf{EA-UVF} & \textbf{EA-UVF1} & \textbf{EA-UVF2} \\
		\hline
		$\#SR$ & \textbf{50/50} & \textbf{50/50} & \textbf{50/50}  & \textbf{50/50} & \textbf{50/50} & 41/50  \\ 
		$M\#G$ & 135.7 & 135.2 & 179.28 & 146.68 & 116.20 & 418.85 \\ 
		$S\#G$ & 90.3  & 77.79 & 100.30 & 90.36  & 56.32  & 299.24 \\
		$A\#P$ & 27.70 & 14.06 & \textit{9.52}   &  & & \\
		$S\#P$ & 18.02 & 7.79  & 5.02   &  & &  \\ 
		$MT$   & 3.26m & 2.48m & 2.63m  &  & &  \\ 
		$ST$   & 3.34m & 1.93m & 1.85m  &   \\	
		$A\_BRSD$ &        &    &    &  & & 0.281    \\	
		\hline
		$U_{1345}^{N}$ & \textbf{NIICh\_5} (PC4) & \textbf{NIICh\_10} (PC4) & \textbf{NIICh\_20} (PC4) & \textbf{EA-UVF} & \textbf{EA-UVF1} & \textbf{EA-UVF2} \\
		\hline
		$\#SR$ & \textbf{50/50}  & \textbf{50/50} & \textbf{50/50}  & \textbf{50/50} & \textbf{50/50}  & 48/50  \\ 
		$M\#G$ & 122.36 & 143.40 & 160.18 & 107.04 & 119.66 & 276.42  \\ 
		$S\#G$ & 64.40  & 65.46 & 77.33  & 62.66 & 66.93 & 212.29  \\
		$A\#P$ & 25.14  & 14.90 & \textit{8.60}   & & &  \\
		$S\#P$ & 12.83  & 6.50  & 3.85   & & &  \\ 
		$MT$   & 2.61m  & 2.66m & 2.68m  & & &  \\ 
		$ST$   & 1.95m  & 1.77m & 1.45m  & & &  \\		
		$A\_BRSD$ &        &    &    &  & & 0.213    \\	
		\hline
		$U_{2345}^{N}$ & \textbf{NIICh\_5} (PC2) & \textbf{NIICh\_10} (PC2) & \textbf{NIICh\_20} (PC2) & \textbf{EA-UVF} & \textbf{EA-UVF1} & \textbf{EA-UVF2} \\
		\hline
		$\#SR$ & \textbf{50/50}  & \textbf{50/50} & \textbf{50/50}  & \textbf{50/50} & \textbf{50/50} & 44/50  \\  
		$M\#G$ & 134.58 & 151.80 & 162.50  & 113.06 & 124.10 & 303.77 \\ 
		$S\#G$ & 79.73  & 60.40 & 67.49  & 82.24 & 85.86 & 229.71  \\
		$A\#P$ & 27.56  & 15.70 & \textit{8.72}   & & &  \\
		$S\#P$ & 15.92  & 5.97  & 3.36   & & &  \\ 
		$MT$   & 2.48m  & 2.03m & 1.66m  & & &   \\ 
		$ST$   & 2.27m  & 1.07m & 45.84s & & &  \\		
		$A\_BRSD$ &        &    &    &  & & 0.257    \\	
		\end{tabular}		
		}		
\end{table}
}

{
\begin{table}[!h]
	\centering
\caption{Results for functions $U^{D}$ and $U_{v}^D$ considering $p=4$. All simulations have been performed with four {different} PCs which characteristics and labels are the following: (PC1) intel core i7 3.6GHz; (PC2) intel core i5 2.5GHz; (PC3) intel core i7 2.7GHz; (PC4) intel core i7 1.9GHz.\label{D4Table}}
\resizebox{0.7\textwidth}{!}{
		\begin{tabular}{ccccccc}
		\hline
    $U^{D}$ & \textbf{NIICh\_5} (PC3) & \textbf{NIICh\_10} (PC3) & \textbf{NIICh\_20} (PC1) & \textbf{EA-UVF} & \textbf{EA-UVF1} & \textbf{EA-UVF2}\\
		\hline
		$\#SR$ &  \textbf{50/50}  &  \textbf{50/50}  & \textbf{50/50}  & 46/50  & \textbf{50/50} & 27/50 \\
		$M\#G$ &  223.16 &  216.64 & 211.72 & 275.72 & 221.40 & 340.52 \\ 
		$S\#G$ &  163.75 &  145.92 & 178.78 & 192.93 & 135.63 & 268.00 \\
		$A\#P$ &  45.30  &  22.22  & \textit{11.14}  &  & & \\
		$S\#P$ &  32.76  &  14.62  & 8.91   &  & & \\ 
		$MT$   &  1.86h  &  1.08h  & 20m    &  & &   \\ 
		$ST$   &  3.41h  &  2.36h  & 55.24m &  & &   \\ 
		$A\_BRSD$ &         &         &        & 0.217 & & 0.091  \\
    \hline
		$U_{1234}^{D}$ & \textbf{NIICh\_5} (PC3) & \textbf{NIICh\_10} (PC1) & \textbf{NIICh\_20} (PC2) & \textbf{EA-UVF} & \textbf{EA-UVF1} & \textbf{EA-UVF2}\\
		\hline
		$\#SR$ &  \textbf{50/50}  & \textbf{50/50}  & \textbf{50/50}  & 48/50 & \textbf{50/50} & 49/50  \\
		$M\#G$ &  133.12 & 175.44 & 228.20 & 176.85 & 156.94 & 214.57 \\ 
		$S\#G$ &  86.49  & 136.76 & 229.41 & 115.39 & 123.64 & 182.38 \\
		$A\#P$ &  27.22  & 18.06  & \textit{12.00}  & & &  \\
		$S\#P$ &  17.32  & 13.68  & 11.45  & & &  \\ 
		$MT$   &  19.81m & 27.68m & 11.29m & & &    \\ 
		$ST$   &  37.04m & 1.33h  & 49.2m  & & &   \\ 
		$A\_BRSD$ &         &         &        & 0.408 & & 0.258   \\
		\hline
		$U_{1235}^{D}$ & \textbf{NIICh\_5} (PC3) & \textbf{NIICh\_10} (PC2) & \textbf{NIICh\_20} (PC2) & \textbf{EA-UVF} & \textbf{EA-UVF1} & \textbf{EA-UVF2}\\
		\hline
		$\#SR$ &  \textbf{50/50}  & \textbf{50/50}  & \textbf{50/50}  & 46/50  & \textbf{50/50} & 22/50 \\
		$M\#G$ &  231.38 & 217.38 & 210.72 & 274.98 & 219.70 & 285.73 \\ 
		$S\#G$ &  178.15 & 147.69 & 178.45 & 183.33 & 144.37 & 277.69 \\
		$A\#P$ &  46.94  & 22.30  & \textit{11.10}  &  & & \\
		$S\#P$ &  35.63  & 14.81  & 8.90   &  & & \\ 
		$MT$   &  14.84h & 4.76h  & 16.49m &  & &   \\ 
		$ST$   &  33.23h & 8.01h  & 49.31m &  & &   \\ 
		$A\_BRSD$ &         &        &        & 0.217 & & 0.074   \\
		\hline
		$U_{1245}^{D}$ & \textbf{NIICh\_5} (PC1) & \textbf{NIICh\_10} (PC1) & \textbf{NIICh\_20} (PC2) & \textbf{EA-UVF} & \textbf{EA-UVF1} & \textbf{EA-UVF2}\\
		\hline
		$\#SR$ &  \textbf{50/50}  & \textbf{50/50}  & \textbf{50/50}  & 21/50 & 39/50 & 47/50   \\
		$M\#G$ &  189.96 & 217.16 & 194.36 & 246.52 & 204.44 & 281.62 \\ 
		$S\#G$ &  156.31 & 226.02 & 151.18 & 254.00 & 178.24  & 248.93 \\
		$A\#P$ &  38.66  & 22.24  & \textit{10.24}  &  & & \\
		$S\#P$ &  31.24  & 22.57  & 7.61   &  & & \\ 
		$MT$   &  13.41h & 2.08h  & 22.73m &  & &  \\ 
		$ST$   &  29.39h & 5.73h  & 1.03h  &  & &   \\ 
		$A\_BRSD$ &         &        &        & 0.079 & 0.078 & 0.095 \\
		\hline
		$U_{1345}^{D}$ & \textbf{NIICh\_5} (PC2) & \textbf{NIICh\_10} (PC2) & \textbf{NIICh\_20} (PC2) & \textbf{EA-UVF} & \textbf{EA-UVF1} & \textbf{EA-UVF2}\\
		\hline
		$\#SR$ &  \textbf{50/50}  & \textbf{50/50}  & \textbf{50/50}  & 46/50 & \textbf{50/50} & 27/50   \\
		$M\#G$ &  223.16 & 216.84 & 210.72 & 275.72 & 220.56 & 340.52 \\ 
		$S\#G$ &  163.75 & 146.07 & 178.45 & 192.93 & 136.10 & 268.00 \\
		$A\#P$ &  45.30  & 22.24  & \textit{11.10}  &  & &  \\
		$S\#P$ &  32.76  & 14.63  & 8.90   &  & &  \\ 
		$MT$   &  2.78h  & 1.38h  & 14.45m &  & &  \\ 
		$ST$   &  5.19h  & 3.1h   & 44.94m &  & &  \\ 
		$A\_BRSD$ &         &        &        & 0.217 & & 0.091   \\
		\hline
		$U_{2345}^{D}$ & \textbf{NIICh\_5} (PC1) & \textbf{NIICh\_10} (PC2) & \textbf{NIICh\_20} (PC1) & \textbf{EA-UVF} & \textbf{EA-UVF1} & \textbf{EA-UVF2}\\
		\hline
		$\#SR$ &  \textbf{50/50}  & \textbf{50/50}  & \textbf{50/50}  & 46/50 & \textbf{50/50} & 27/50  \\
		$M\#G$ &  223.16 & 216.64 & 210.72 & 275.72 & 221.40 & 340.52 \\ 
		$S\#G$ &  163.75 & 145.92 & 178.45 & 192.93 & 135.63 & 268.00 \\
		$A\#P$ &  45.30  & 22.22  & \textit{11.10}  &  & &  \\
		$S\#P$ &  32.76  & 14.62  & 8.90   &  & & \\ 
		$MT$   &  2.38h  & 2.15h  & 16.02m &  & &   \\ 
		$ST$   &  4.33h  & 4.6h   & 47.06m &  & &   \\ 
		$A\_BRSD$ &         &        &        & 0.217 & & 0.091   \\
		\end{tabular}		
		}		
\end{table}
}

{
\begin{table}[!h]
	\centering
\caption{Results for functions $U^{D}$ and $U_{v}^D$ considering $p=5$. All simulations have been performed with four {different} PCs which characteristics and labels are the following: (PC1) intel core i7 3.6GHz; (PC2) intel core i5 2.5GHz; (PC3) intel core i7 2.7GHz; (PC4) intel core i7 1.9GHz.\label{D5Table}}
\resizebox{0.7\textwidth}{!}{
		\begin{tabular}{ccccccc}
		\hline
    $U^{D}$ & \textbf{NIICh\_5} (PC1) & \textbf{NIICh\_10} (PC2) & \textbf{NIICh\_20} (PC3) & \textbf{EA-UVF} & \textbf{EA-UVF1} & \textbf{EA-UVF2} \\
		\hline
		$\#SR$ &  \textbf{45/50}  & 43/50  & 40/50  & 21/50 & 33/50 & 10/50   \\
		$M\#G$ &  433.33 & 365.49 & 455.20 & 214.71 & 508.92 & 337.70 \\ 
		$S\#G$ &  281.33 & 229.93 & 239.48 & 204.39 & 392.73 & 266.71 \\
		$A\#P$ &  87.31  & 37.19  & \textit{23.25}  & & &  \\
		$S\#P$ &  56.25  & 22.95  & \textit{12.00}  & & &  \\ 
		$MT$   &  13.61h & 4.91h  & 1.81h  & & &  \\ 
		$ST$   &  12.58h & 7.27h  & 2.44h  & & &  \\  
		$A\_BRSD$ &  0.038  & 0.093  & 0.104  & 0.103 & 0.101 & 0.175  \\	
		\hline
		$U_{1234}^{D}$ & \textbf{NIICh\_5} (PC2) & \textbf{NIICh\_10} (PC2) & \textbf{NIICh\_20} (PC2) & \textbf{EA-UVF} & \textbf{EA-UVF1} & \textbf{EA-UVF2} \\
		\hline
		$\#SR$ & \textbf{50/50}  & \textbf{50/50}  & \textbf{50/50}  & 46/50 & \textbf{50/50} & 4/50   \\
		$M\#G$ & 185.92 & 208.92 & 304.18 & 199.87 & 169.44 & 472.25 \\ 
		$S\#G$ & 84.31  & 103.11  & 168.53 & 145.30 & 105.37 & 325.81 \\
		$A\#P$ & 37.80  & 21.48  & \textit{15.66}   & & &  \\
		$S\#P$ & 16.86  & 10.34  & 8.41   & & &  \\ 
		$MT$   & 1.74h  & 51.27m & 1.87h  & & &  \\ 
		$ST$   & 1.96h  & 1.69h  & 3.52h  & & &  \\ 
		$A\_BRSD$ &        &        &        & 0.306 & & 0.518   \\
		\hline
		$U_{1235}^{D}$ & \textbf{NIICh\_5} (PC2) & \textbf{NIICh\_10} (PC2) & \textbf{NIICh\_20} (PC2) & \textbf{EA-UVF} & \textbf{EA-UVF1} & \textbf{EA-UVF2} \\
		\hline
		$\#SR$ & \textbf{45/50}  & 44/50  & 40/50  & 21/50  & 30/50 & 8/50  \\
		$M\#G$ & 424.24 & 379.14 & 455.20 & 209.24 & 198.73 & 220.5 \\ 
		$S\#G$ & 278.44 & 247.16 & 239.48 & 191.76 & 123.14 & 142.64 \\
		$A\#P$ & 85.49  & 38.52  & \textit{23.25}  &  & &  \\
		$S\#P$ & 55.70  & 24.70  & 12.00  &  & &  \\
		$MT$   & 9.01h  & 5.52h  & 2.52h  &  & &  \\
		$ST$   & 8.92h  & 7.86h  & 3.41h  &  & &  \\
		$A\_BRSD$ & 0.067  & 0.038  & 0.101  & 0.104 & 0.098 & 0.16   \\
		\hline
		$U_{1245}^{D}$ & \textbf{NIICh\_5} (PC3) & \textbf{NIICh\_10} (PC2) & \textbf{NIICh\_20} (PC1) & \textbf{EA-UVF} & \textbf{EA-UVF1} & \textbf{EA-UVF2} \\
		\hline
		$\#SR$ & \textbf{10/50}  & 5/50   & 9/50   & 5/50  & 4/50 & 5/50     \\
		$M\#G$ & 519.30 & 663.20 & 689.44 & 435.80 & 260.00 & 661.00    \\
		$S\#G$ & 266.92 & 301.91 & 262.90 & 316.05 & 146.81 & 338.02 \\
		$A\#P$ & 104.50 & 66.80  & \textit{34.78}  & & &            \\
		$S\#P$ & 53.30  & 30.33  & 13.12  & & &            \\
		$MT$   & 6.03h  & 10.79h & 14.97h & & &            \\
		$ST$   & 3.84h  & 6.84h  & 11.16h & & &            \\
		$A\_BRSD$ & 0.016  & 0.014  & 0.021  & 0.046 & 0.019 & 0.103 \\
		\hline
		$U_{1345}^{D}$ & \textbf{NIICh\_5} (PC2) & \textbf{NIICh\_10} (PC2) & \textbf{NIICh\_20} (PC2) & \textbf{EA-UVF} & \textbf{EA-UVF1} & \textbf{EA-UVF2} \\
		\hline
		$\#SR$ & \textbf{43/50}  & \textbf{43/50}  & 39/50  & 21/50  & 27/50 & 9/50    \\
		$M\#G$ & 396.88 & 388.47 & 417.90 & 183.52 & 276.37 & 318.33    \\
		$S\#G$ & 258.81 & 229.82 & 234.29 & 115.81 & 203.549 & 293.59    \\
		$A\#P$ & 80.05  & 39.49  & \textit{21.41}  & & &            \\
		$S\#P$ & 51.77  & 22.96  & 11.73  & & &            \\
		$MT$   & 7.58h  & 8.14h  & 2.73h  & & &            \\
		$ST$   & 7.2h   & 11.09h & 3.94h  & & &            \\
		$A\_BRSD$ & 0.045  & 0.086  & 0.093  & 0.108 & 0.103 & 0.16 \\
		\hline
		$U_{2345}^{D}$ & \textbf{NIICh\_5} (PC3) & \textbf{NIICh\_10} (PC3) & \textbf{NIICh\_20} (PC3) & \textbf{EA-UVF} & \textbf{EA-UVF1} & \textbf{EA-UVF2} \\
		\hline
		$\#SR$ & \textbf{45/50}  & 43/50  & 40/50  & 21/50 & 33/50 & 10/50  \\
		$M\#G$ & 433.33 & 365.49 & 455.20 & 214.71 & 255.94 & 337.70 \\
		$S\#G$ & 281.33 & 229.93 & 239.48 & 204.39 & 205.09 & 266.71 \\
		$A\#P$ & 87.31  & 37.19  & \textit{23.25}  & & &  \\
		$S\#P$ & 56.25  & 22.95  & 12.00  & & &  \\
		$MT$   & 5.19h  & 4.61h  & 1.81h  & & &  \\
		$ST$   & 4.58h  & 6.94h  & 2.44h  & & &  \\
		$A\_BRSD$ & 0.038  & 0.075  & 0.103  & 0.103 & 0.101 & 0.175 \\
		\end{tabular}		
		}		
\end{table}
}

{
In Tables \ref{N4Table}-\ref{D5Table} we reported the results of the application of the three versions of NEMO-II-Ch as well as those obtained by the three algorithms knowing the user's true value function. We have considered the twelve different user's true value functions defined in the previous section and the cases $p=4$ and $p=5$ for the number of best locations to be discovered. 
}

{
\noindent In the following, by $\left(U,p\right)$ we denote the case in which the user's true value function is $U$ and the number of best locations is $p$. The following can be observed: \\
\begin{itemize}
\item $\left(U^{N},4\right)$ and $\left(U_{v}^{N},4\right)$:
\begin{itemize}
\item \textit{Convergence:} the three variants of NEMO-II-Ch as well as EA-UVF and EA-UVF1 are always able to find the best solution in the 50 runs. This is not the case for EA-UVF2 that, with respect to $U^N$ is not able to converge in one of the 50 runs, while, with respect to $U^N_{1245}$, quite surprisingly, it is able to find the best subset of 4 locations only in 5 of the 50 runs; 
\item \textit{Convergence speed:} As can be observed from the data in Table \ref{N4Table}, apart from the $U^{N}$ case in which the EA-UVF1 converges more quickly (in terms of number of generations necessary to find $P_b$) than all the other algorithms, the EA-UVF is the quickest among the considered algorithms. As to the comparison between the three NEMO-II-Ch variants, in average, NIICh\_5 converges more quickly than NIICh\_10 in four of the six considered cases, while NIICh\_20 is always the slowest. However, as already observed before, the number of pairwise comparisons asked from EA-UVF and EA-UVF1 is tremendously higher than the one involved in whichever NEMO-II-Ch version. For this reason, it is more meaningful giving a more in depth analysis of the NEMO-II-Ch variants to understand if and how the number of times the user is queried with a pairwise comparison affects the convergence speed of the algorithm. It can be observed that the lowest number of pairwise comparison is asked in correspondence of NIICh\_20, followed by NIICh\_10 and, then, by NIICh\_5 (see values in italics). This means that not only NIICh\_20 is efficient in finding $P_b$ but it is able to find it asking very few pairwise comparisons to the user; 
\item \textit{Distance from $P_b$:} Considering EA-UVF2 and assuming that the best solution in the final population is the optimal one, the user makes an error, in average, of the 40.6\% in the $U^{N}$ case, and of the 23.3\% in the $U_{1245}^{N}$ one;
\end{itemize}
\item $\left(U^{N},5\right)$ and $\left(U_{v}^{N},5\right)$:
\begin{itemize}
\item \textit{Convergence:} The three variants on NEMO-II-Ch are able to find $P_b$ in all considered runs for all test problems apart from the case $\left(U^{N}_{1235}\right)$ in which NIICh\_10 and NIICh\_20 are not always able to find $P_b$. In particular, NIICh\_10 does not find the best subset of five locations in one of the 50 runs, while NIICh\_20 does not find the same subset of best locations in 3 out of the 50 runs. \\
As to the three algorithms knowing the user's true value functions, EA-UVF1 is always able to find the best subset of five locations, while this is not true for the other two. In particular, EA-UVF does not find $P_b$ in five of the fifty runs in the $U^N_{1235}$ case, while EA-UVF2 has its best performances when $U^{N}$ is considered (49/50) and its worst one in the case $U_{1235}^{N}$ is the user's true value function (11/50). This suggests that using the user's true value function to assign a probability to become parent of the next generation is worse than using the same function to rank the solutions belonging to the same front; 
\item \textit{Convergence speed:} As in the $p=4$ cases, it results that NIICh\_20 is the quickest among the three NEMO-II-Ch versions to reach $P_b$ since, it asks to the user to provide almost half of the pairwise comparisons asked by NIICh\_10 and almost one third of the pairwise comparisons asked by NIICh\_5. \\
Regarding EA-UVF and its two variants, once again EA-UVF2 is the worst among them. Moreover, we would like to underline that the number of generations necessary to get $P_b$ is lower for EA-UVF1 than for EA-UVF. In particular, it is meaningful observing that the number of pairwise comparisons asked to the user by EA-UVF1 is not greater than the number of times the user is queried with a pairwise comparison in the EA-UVF. In fact, the application of the EA-UVF1 implies the same number of pairwise comparisons of EA-UVF only in case all solutions are non-dominated and, therefore, they are in one non-dominated front only. This suggests once again that a parsimonious preference information is beneficial for the convergence of the algorithms to $P_b$;
\item \textit{Distance from $P_b$:} In the $U_{1235}^{N}$ case, in average, the error done in assuming that the best solution in the last population is the optimal one is almost 7\% for NIICh\_10, NIICh\_20 and EA-UVF, while it is 21.3\% for the EA-UVF2. An higher error is also done by EA-UVF2 in the $U_{1234}^{N}$, $U_{1245}^{N}$ and $U_{2345}^{N}$ cases.
\end{itemize}
\item $\left(U^{D},4\right)$ and $\left(U_{v}^{D},4\right)$:
\begin{itemize}
\item \textit{Convergence:} The three variants of NEMO-II-Ch are always able to find $P_b$ in all considered runs. This is not the case for the three algorithms knowing the user's true value function. In particular, EA-UVF1 finds in all 50 runs the best subset of four locations for all user's true value functions apart from $U^{D}_{1245}$ in which it finds $P_b$ in 39 of the 50 runs; the EA-UVF never finds the best subset of locations in all runs. The same holds for EA-UVF2 that in the $U^{D}_{1234}$ and $U^{D}_{1245}$ cases finds $P_b$ 49 and 47 times, respectively. Considering all other user's true value functions, it is able to find the best subset of four locations more or less half of the times;
\item \textit{Convergence speed:} NIICh\_20 is confirmed as the best among the three variants of the NEMO-II-Ch since it asks a lower number of pairwise comparisons than the other two maintaining always the best possible convergence since, as observed in the previous item, it is always able to find $P_b$. Comparing NIICh\_10 and NIICh\_5, the first is better than the second in terms of number of pairwise comparisons asked to the DM;
\item \textit{Distance from $P_b$:} Assuming that the best solution in the last population is optimal, one makes an error ranging from 7.9\% to 40.8\% considering the EA-UVF, from 7.4\% to 25.8\% considering the EA-UVF2 and of the 7.8\% considering the EA-UVF1. 
\end{itemize}
\item $\left(U^{D},5\right)$ and $\left(U_{v}^{D},5\right)$:
\begin{itemize}
\item \textit{Convergence:} For all considered cases, one of the three variants of NEMO-II-Ch finds $P_b$ more often than the algorithms based on the knowledge of the user's true value function. Even more, in all cases the worst among the three NEMO-II-Ch variants performs at least as well as all three algorithms knowing the user's true value function in terms of number of runs in which it converges to $P_b$; 
\item \textit{Convergence speed:} Looking at the average number of pairwise comparisons asked to the user, once more we have the confirmation that NIICh\_20 is the best among the three variants of NEMO-II-Ch since it finds $P_{b}$ asking less pairwise comparisons than NIICh\_5 and NIICh\_10. However, differently from the previous cases, the doubt is now related to the fact that NIICh\_20 is not able to find the best subset of locations as frequently as NIICh\_5 and NIICh\_10 and, therefore, it could be better to ask more pairwise comparisons to increase the probability to converge to the best solution. 
\item \textit{Distance from $P_b$:} Comparing the three versions of NEMO-II-Ch one can observe that, apart from $U^{D}_{1235}$ and $U^{D}_{1245}$ cases, NIICh\_5 presents the best $A\_BRSD$. In particular, the maximum average error is equal to 6.7\% for NIICh\_5, while it is 9.3\% for NIICh\_10 and even 10.4\% for NIICh\_20. The situation is even worse for the three algorithms based on the full knowledge of the user's true value function since, apart from the $U_{1245}^{D}$ case in which the average error done assuming as optimal solution the best solution in the final population is 1.9\% considering EA-UVF1 and 4.6\% considering EA-UVF, in all the other cases, this average error is at least equal to 9.8\% with a pick of 51.8\% done by EA-UVF2 in the $U^{D}_{1234}$ case. This means that, in the case in which the EA-UVF algorithm and the other two variants are not able to find $P_b$, they are very far from the area of the Pareto front most interesting with respect to the user's preferences. 
\end{itemize}
\end{itemize}
}

{To evaluate the significance of the data provided above we performed the Mann-Whitney $U$ test with $5\%$ significance level \citep{HollanderWolfeChicken2013} to two different indicators: 
\begin{enumerate}
\item considering $BRSD$ of each of the six algorithms in each of the 50 considered runs, 
\item considering the number of pairwise comparisons asked to the user in each run for algorithms NIICh\_5, NIICh\_10 and NIICh\_20. 
\end{enumerate}
}

{Regarding the $BRSD$, we performed the test only for problems where at least one algorithm did not converge in at least one of the 50 runs. Indeed, if all methods had converged to the optimal solution in all runs, the $BRSD$ would be always equal to 0 and, consequently, the comparison between the algorithms would be absolutely meaningless. \\
Regarding the number of pairwise comparisons asked to the user, we performed the test on the NIICh\_5, NIICh\_10 and NIICh\_20 only, since the number of pairwise comparisons asked to the user in EA-UVF, EA-UVF1 and EA-UVF2 is only virtual due to the unrealistic applicability of the algorithms. In particular, in the case in which the algorithm did not converge to the optimal solution, for that run, we considered the maximum number of pairwise comparisons asked to the user being 200 for NIICh\_5, 100 for NIICh\_10 and 50 for NIICh\_20 since each of them asks one pairwise comparison every 5, 10 and 20 generations, respectively, and the maximum number of admitted generations is 1,000. }

{In the supplementary material we included the results of the two tests. For brevity, we report here just the tables for the $\left(U^{D},5\right)$ case obtained performing the Mann-Whitney $U$ test with $5\%$ significance level on the $BRSD$ (Table \ref{MWTest5D}) and on the number of pairwise comparisons asked to the user (Table \ref{MWTestAverPref5D}). In both tables, we give the $p$-value together with the difference between the $A\_BRSD$ of each ordered pair of algorithms in Table \ref{MWTest5D} and the difference between $A\#P$ of each ordered pair of algorithms in Table \ref{MWTestAverPref5D}. Red values represent significant values considering the performed test. }

\begin{table}[!h]
	\centering
\caption{Mann-Whitney $U$ test with $5\%$ significance level performed on $BRSD$ for the $\left(U^D,5\right)$. In the table the $p$-value is provided as well as the difference between the $A\_BRSD$ of each ordered pair of algorithms. In red the significant values.\label{MWTest5D}}
\resizebox{0.8\textwidth}{!}{
{\renewcommand\arraystretch{1.3} 
		\begin{tabular}{ccccccc}
		\hline
    $U^{D}$ & \textbf{NIICh\_5} & \textbf{NIICh\_10} & \textbf{NIICh\_20} & \textbf{EA-UVF} & \textbf{EA-UVF1} & \textbf{EA-UVF2} \\
		\hline
		\textbf{NIICh\_5}  & & $\substack{0.4598 \\ (0.0038-0.013)}$ & $\substack{0.1137 \\ (0.0038-0.0209)}$ & $\substack{\textcolor{red}{5.28\cdot 10^{-8}} \\ (0.0038-0.0595)}$ & $\substack{\textcolor{red}{0.0013} \\ (0.0038-0.0345)}$ & $\substack{\textcolor{red}{4.10\cdot 10^{-13}} \\ (0.0038-0.1398)}$ \\
		\textbf{NIICh\_10} &   &  & $\substack{0.4355 \\ (0.0130-0.0209)}$ & $\substack{\textcolor{red}{9.59\cdot 10^{-6}} \\ (0.013-0.0595)}$ & $\substack{\textcolor{red}{0.0220} \\ (0.0130-0.0345)}$ & $\substack{\textcolor{red}{3.31\cdot 10^{-11}} \\ (0.013-0.1398)}$ \\
		\textbf{NIICh\_20} &   &  &  & $\substack{\textcolor{red}{0.0002} \\ (0.0209-0.0595)}$ & $\substack{{0.1517} \\ (0.0209-0.0345)}$ & $\substack{\textcolor{red}{6.42\cdot 10^{-10}} \\ (0.0209-0.1398)}$ \\
		\textbf{EA-UVF}    &   &  &  &  & $\substack{\textcolor{red}{0.0120} \\ (0.0595-0.0345)}$ & $\substack{\textcolor{red}{0.0001} \\ (0.0595-0.1398)}$ \\
		\textbf{EA-UVF1}   &   &  &  &  &  &  $\substack{\textcolor{red}{1.36\cdot 10^{-7}} \\ (0.0345-0.1398)}$ \\
		\hline
		\end{tabular}		}
		}		
\end{table}

\begin{table}[!h]
	\centering
\caption{Mann-Whitney $U$ test with $5\%$ significance level on the number of pairwise comparisons asked to the user in algorithms NIICh\_5, NIICh\_10 and NIICh\_20 for the $\left(U^D,5\right)$ case. In the table the $p$-value is provided as well as the difference between $A\#P$ of each ordered pair of algorithms. In red the significant values.\label{MWTestAverPref5D}}
\resizebox{0.5\textwidth}{!}{
{\renewcommand\arraystretch{1.3} 
		\begin{tabular}{cccc}
\hline
$U^{D}$ & \textbf{NIICh\_5} & \textbf{NIICh\_10} & \textbf{NIICh\_20} \\
\hline
\textbf{NIICh\_5}  & & $\substack{ \textcolor{red}{1.82\cdot 10^{-5}} \\ (98.68-46.12)}$ & $\substack{ \textcolor{red}{2.12\cdot 10^{-9}} \\ (98.68-28.8) }$ \\
\textbf{NIICh\_10} & & & $\substack{ \textcolor{red}{0.0105} \\ (46.12-28.8)}$ \\
\hline
\end{tabular}		}
		}		
\end{table}

{In Table \ref{MWTest5D} one can observe that the difference in the $BRSD$ between the NEMO variants is not significant, while the difference between the $BRSD$ of each NEMO variant and each of the algorithms based on the knowledge of the user's true value function is significant apart from the comparison between NIICh\_20 and EA-UVF1 for which the difference between their $BRSD$ is not significant for the Mann-Whitney $U$ test. This means that, on the one hand, the NEMO-II-Ch variants can be considered equivalent, while each of them is better than the three algorithms knowing the user's true value function. On the other hand, one can conclude that EA-UVF1 is better than EA-UVF that, in turn, is better than EA-UVF2.} 

{Going to the data in Table \ref{MWTestAverPref5D}, one can see that the difference between the distributions of the number of pairwise comparisons asked to the user in each pair of NEMO variants is significant. This means that the number of pairwise comparisons asked to the user by NIICh\_20 to converge to the optimal solution is retained significantly smaller than the one involved in NIICh\_10 and NIICh\_5 and, consequently, with respect to the required preference information, NIICh\_20 is better than NIICh\_10 that, in turn, is better than NIICh\_5.}

{Similar conclusions can be gathered looking at all the other tables included in the supplementary material. Once again they confirm that the difference in the $BRSD$ between the NEMO variants and the algorithm based on the user's true value function is considered significant and that with respect to the three NEMO variants, the difference between the number of pairwise comparisons asked to the user from each algorithm is significant. The last fact proves that asking to the user a lower number of information does not affect, in general, the capacity of NEMO-II-Ch to converge to the optimal solution. }

\section{Discussion}\label{DiscussionSection}
In this paper we faced Multiobjective Combinatorial Optimization (MOCO) problems by using Interactive Evolutionary Multiobjective Optimization (IEMO) methods. In particular, we applied an IEMO method, namely NEMO-II-Ch, to Facility Location Problems (FLPs). When facility location is considered from a multiobjective perspective, that is the different facility location options are evaluated simultaneously by several conflicting aspects, the problem becomes quite difficult and the user is presented with a huge number of Pareto optimal solutions from which he is asked to choose the best w.r.t. his preferences. Anyway, this set can be composed of thousands and even millions of different possible facility locations options and, therefore, no user can cope with such a type of problem and in a reasonable time. Moreover, the user has not clear defined preferences at the beginning of the search process, so that, the already difficult multiobjective optimization problem, is {coupled} with a not less complex problem of preference elicitation. In fact, what we are proposing can be considered {as} the first structured methodology in FLPs to search optimal solutions taking into account preferences of the user.\\
For these reasons, we proposed to apply NEMO-II-Ch that is a state of the art algorithm permitting to conjugate the efficiency of evolutionary multiobjective optimization methods with the parsimonious user's preferences elicitation of most advanced multiple criteria decision aiding models, with the aim of searching the best solution w.r.t. the preferences of the user in the part of the Pareto front most appealing for him. In this way, focusing on a single region of the Pareto front, the algorithm reaches the best option in a limited number of generations avoiding to lose time in searching solutions in parts of the Pareto front not really interesting for the user.

To prove the efficiency of the method to this setting, we considered a classical FLP very well-known in literature \cite{drezner2006multi} based on the most typical objective functions adopted in the domain. We performed different simulations running NEMO-II-Ch and comparing its performances with those of other three algorithms, namely EA-UVF, EA-UVF1 and EA-UVF2, based on the knowledge of the user's true value function that is, instead, unknown to NEMO-II-Ch. \\
In the comparison we tested twelve different types of users' value functions and two different values for the number of facilities $p$ that need to be located ($p=4$ and $p=5$). Moreover, to investigate how the number of comparisons asked to the user influences the convergence of the algorithm, we considered three different versions of the NEMO-II-Ch method, namely NIICh\_5, NIICh\_10 and NIICh\_20, asking the user to compare one pair of non-dominated solutions every 5, 10 and 20 generations, respectively.


{The results obtained should be read as an answer to the question: ``\textit{is there any methodological tool to handle real world multiobjective facility problems}"? The considered problem is very complex for the following reasons:
\begin{enumerate}
\item There is a plurality of objectives to be optimized,
\item Some of these objectives are quite complex in itself (this is, in particular, the case of $f_5(P)$ \cite{drezner2006multi}), 
\item The preferences of the user have to be considered,
\item The preference information has to be collected maintaining tolerable the cognitive burden for the DM,
\item The computation time should be acceptable for real world operational applications.
\end{enumerate}
The data of the conducted simulations allow to conclude that the proposed methodology answers positively the research question. Indeed, NEMO-II-Ch is able to find the optimal solution in most of the considered FLPs taking into account the preferences of the user, without requiring too much preference information and involving a computational time definitely admissible. The performances of the ``control" procedures having full preference information, that is, EA-UVF, EA-UVF1 and EA-UVF2, can be used to measure the task complexity in the sense that the worse their performances, the more complex the task. In particular, let us evaluate the complexity of the problem taking into account EA-UVF1, having the best performances between the three algorithms knowing the user's true value function. Considering as object of the performance the number of runs in which the optimal solution was obtained, that is indicator \#SR, we can see that almost always NIICh\_20 obtains performances at least as good EA-UVF1 and sometimes even better. More precisely NEMO-II-Ch is obtaining better performances for the cases 
\begin{itemize}
\item $\left(U^D_{1245},4\right)$ (NIICh\_20 found the solution in all the 50 runs, while EA-UVF1 in 39 runs), 
\item $\left(U^D,5\right)$ (NIICh\_20 found the solution in 40 runs, while EA-UVF1 in 33 runs),
\item $\left(U^D_{1235},5\right)$ (NIICh\_20 found the solution in 40 runs, while EA-UVF1 in 30 runs), 
\item $\left(U^D_{1245},5\right)$ (NIICh\_20 found the solution in 9 runs, while EA-UVF1 in 4 runs),
\item $\left(U^D_{1345},5\right)$ (NIICh\_20 found the solution in 39 runs, while EA-UVF1 in 27 runs),
\item $\left(U^D_{2345},5\right)$, (NIICh\_20 found the solution in 40 runs, while EA-UVF1 in 33 runs).
\end{itemize}
Instead, there is only one case in which EA-UVF1 is performing better than NIICh\_{20} with respect to the number of runs in which the optimal solution is found, that is $\left(U^N_{1235},5\right)$ (NIICh\_{20} found the optimal solution in 47 runs, while EA-UVF1 in all 50 runs). In all other cases, both algorithms where able to find the optimal solution in all runs.  }

{Considering the number of pairwise comparison requested by NEMO-II-Ch we have to conclude that it is definitely acceptable. Indeed,in terms of comparison the pairwise comparisons requested by one of the most well-known and most adopted MCDA method, that is AHP \cite{Saaty1980}. Let us consider the didactic example presented in \cite{Saaty1977} in which three schools (alternatives) are evaluated with respect to six different aspects (criteria). In terms of FLPs, it would be a really easy problem that will concern the selection of a single facility among three potential locations to optimize six different objectives. Since the DM must provide a pairwise comparison in terms of a qualitative judgment on a nine point scale for each non-ordered pair of criteria and a comparison for each non-ordered pair of alternatives with respect to each criterion, the decision maker has to provide $\binom{6}{2}+6\binom{3}{2}=15+6\cdot 3=33$ pairwise comparisons in total. This means that in a didactic example of, probably, the most adopted MCDA method \citep{VaidyaKumar2006}, the DM is asked to give 33 pairwise comparisons. Looking again at the performances of NIICh\_20, one can see that with a single exception, in all our cases the algorithm was able to find the best solution with a number of pairwise comparison much smaller than 33. Observe also that very often the required average number of pairwise comparison asked to the user by NIICh\_20 is lower than 15 (in 17 out of 24 considered cases). In addition, observe that while the pairwise comparisons of AHP require to give an evaluation on a nine point scale, the pairwise comparisons considered in NEMO-II-Ch require simply to say which solution is the preferred among two. To have a more fair comparison between the judgments required by AHP and the information required by NEMO-II-Ch, consider that for each pair of items $\alpha$ and $\beta$ being alternatives ($\alpha,\beta \in A$) or criteria ($\alpha,\beta \in G$)  AHP requires, in fact, two comparisons: the first related to which one between $\alpha$ and $\beta$ has the greatest priority and the second, expressed on the nine point scale,  related to how much greater is the priority of the item with the greatest priority with respect to the other. In general, it seems reasonable that the second comparison of AHP (the one on the nine point scale) is more demanding than the pairwise comparison of NEMO-II-Ch related to which solution is the preferred among two. Consequently, to each pairwise comparison asked from AHP on a pair of non-ordered items one should assign a cognitive burden at least double with respect to the pairwise comparison required by NEMO-II-Ch. In conclusion, we can say that, in average, NEMO-II-Ch can handle a quite challenging problem with a complexity comparable to that one of the most demanding real world problems asking the user a cognitive burden much smaller than the one required by the most adopted MCDA method in a very didactic example.\\
Coming to the computational time, even considering the case taking more time, NIICh\_20 is almost always (apart from one case only) achieving the optimal solution, in average, in less than three hours and, very often, in less than one hour (quite frequently in the $U^{N}$ and $U^{N}_{v}$ cases in some minutes). This seems a very reasonable running time for a so complex problem. Observe also that our simulations were performed with a non dedicated programming language and with computer machine commonly available on the market as will be further underlined below. }

{Beyond the specific interest for the multiobjective facility location problems, the results we obtained are relevant also from the general point of view of the multiobjective optimization algorithms. In fact, the procedure that has been proposed can be seen as a parsimonious exploration of the space of solutions and of the DM's preferences.  The parsimony of the multiobjective optimization procedure we have applied can be decomposed in two components:
\begin{itemize}
\item	a component related to the optimization procedure: it is based on the evaluations of combinations of most promising solutions maintaining a certain level of diversification typical of the evolutionary algorithms, 
\item	a component related to the preference learning procedure: it is based on a ``dynamical" induction of the DM's utility function on the basis of few preference comparisons, typical of the ordinal regression approach \citep{jacquet1982assessing} that is properly applied in an ``incremental" version adding time by time preferences related to new solutions discovered by the optimization algorithm.
\end{itemize}
The parsimony of the optimization algorithm permits to select the most promising directions in the search of the optimal solution avoiding to be trapped in some local optimum, while the parsimony of the preference learning algorithm permits to escape from a possible overfitting originated by an excess of preference information that could bring to an inappropriate generalization at global level of preferences that hold only at local level. With respect to the parsimony of the preference learning algorithm, let us point out that the obtained results have an autonomous interest. Indeed, the comparisons of the results obtained by NIICh\_20 with the three algorithms based on a complete knowledge of the DM's preferences, proves that a quite limited use of preference information gives better results than the use of the complete preference information. In this sense, it is particularly meaningful that, in general, the best performances are obtained by EA-UVF1. Indeed, while EA-UVF, based only on the complete and perfect preference information, is operating as a single objective evolutionary algorithm without taking into account diversity, and EA-UVF2, beyond considering the multiobjective nature of the problem through dominance front ranking, is based on the maintenance of the diversity by means of crowding distance, EA-UVF1 is obtained as a ``mutiobjectivization" of EA-UVF building dominance front ranking without any consideration of the diversification. This means that, according to a concordant literature in the evolutionary optimization domain \cite{KnowlesWatsonCorne2001,SeguraEtAl2016}, the multiobjecivization of the optimization problem is beneficial, while the diversification, in general, does not give any contribution to improve the performances of the algorithm. In particular, taking into consideration the number of runs in which the best solution was discovered, observe that, on the one hand, EA-UVF2 is performing better than EA-UVF only in two cases ($\left(U^D_{1234},4\right)$ and $\left(U^D_{1245},4\right)$), while EA-UVF is performing better than EA-UVF2 in 17 cases and they have the same performances on the remaining 5 cases. On the other hand, EA-UVF1 is performing better than EA-UVF in 12 cases, while EA-UVF is performing better than EA-UVF1 only in one case ($\left(U^D_{1245},5\right)$ with 5 successful runs for EA-UVF and 4 successful runs for EA-UVF1). Observe that if the multiobjectivization of the EA-UVF1 algorithm permits to improve the performance of EA-UVF, however, it is not enough to attain the same efficiency of NIICh\_20. In fact, taking into consideration the number of runs in which the optimal solution was discovered, NIICh\_20 is able to obtain better results than EA-UVF1 in 6 cases, while EA-UVF1 is able to perform better than NIICh\_20 in one case only. We believe that this can be interpreted in the sense that the parsimony in the required preference information of  NIICh\_20 permits to obtain better performances of an algorithm using the whole preference information as EA-UVF1. In conclusion, the results we obtained on the multiobjective facility location problem seems to suggest that in very complex combinatorial optimization problems a smart approach based on an evolutionary algorithm and a limited elicitation of the DM's preference information can be an appropriate approach. Of course, this hypothesis needs to be tested on other multiobjective combinatorial problems and, more in general, on other complex multiobjective problems (not necessarily combinatorial), to obtain a more precise and definitive confirmation.}
 
\section{Conclusions}
{We considered a very complex problem resulting from the combination of two other complex problems already quite challenging in themselves. The combination of the two problems highly exacerbates the difficulty. The two problems are the facility location problem and
the search of optimal solutions in multiobjective decision problems taking into account the user's preferences. In this perspective, the research question of the paper is: ``is it possible to give an adequate answer, especially taking into account real world applications, to the so complex problem resulting from the combination of the above-mentioned problems?" Technically the answer to the problem is obtained from the application of a state-of-the-art multiobjective optimization procedure to the standard formulation of multiobjective facility location problem. The contribution of the paper is in handling the question and in providing a surprisingly very positive answer: the two complex problems can be solved together with a reasonable cognitive burden (comparable and even
smaller than the cognitive burden required from didactic examples of the most adopted MCDA methods) and with reasonable computational times (especially considering the use of non-specialised programming languages and the computation on common laptops daily used). } Beyond the application of the presented methodology to other complex multiobjective  combinatorial optimization problems in order to collect further evidence on its effectiveness and reliability in so complex decision problems, the following possible avenues of research can be underlined: 
\begin{itemize}
\item research should be addressed on studying how often the user should be asked to provide preference information to speed the convergence of the algorithm and how techniques investigating which solutions should be presented to the user to maximize the learning capabilities of the algorithm \cite{bcggCOR,CiomekKadzinskiTervonen:2017} could improve the same convergence;
\item to make applicable to big size real world problems, a better implementation of NEMO-II-Ch should be provided. Indeed, analyzing in detail the computational time necessary to run the algorithm, it is evident that almost 93\% of the time is taken by the execution of the Nelder-Mead method. Of course, implementing other methods to solve non-linear optimization problems could speed the algorithm and, therefore, making it more applicable in practice,
\item on the basis of the good results obtained by NEMO-II-Ch applied to location problems, we think that it could be interesting applying it to other classical combinatorial optimization problems that can be formulated in a multiobjective perspective such as the one presented in \cite{doerner2009multi} and \cite{Hamacher:2002}.
\end{itemize}

\section*{Acknowledgements}
\noindent The authors are  grateful to Professor Tammy Drezner for making available the data concerning the real world problem. The second and the third authors wish to acknowledge the support of the Ministero dell'Istruzione, dell'Universit\'{a} e della Ricerca (MIUR) - PRIN 2017, project "Multiple Criteria Decision Analysis and Multiple Criteria Decision Theory", grant 2017CY2NCA. Salvatore Corrente wishes to acknowledge also the support of the Starting Grant 2020 from the University of Catania.


\bibliographystyle{plainnat}
\bibliography{Full_bibliography}

\begin{thebibliography}{84}
\providecommand{\natexlab}[1]{#1}
\providecommand{\url}[1]{\texttt{#1}}
\expandafter\ifx\csname urlstyle\endcsname\relax
  \providecommand{\doi}[1]{doi: #1}\else
  \providecommand{\doi}{doi: \begingroup \urlstyle{rm}\Url}\fi

\bibitem[Alcada-Almeida et~al.(2009)Alcada-Almeida, Coutinho-Rodrigues, and
  Current]{ALCADAALMEIDA2009111}
L.~Alcada-Almeida, J.~Coutinho-Rodrigues, and J.~Current.
\newblock {A multiobjectrive modeling approach to Locating incinerators}.
\newblock \emph{Socio-Economic Planning Sciences}, 43\penalty0 (2):\penalty0
  111--120, 2009.

\bibitem[Alves and Cl{\'\i}maco(2007)]{alves2007review}
M.J. Alves and J.~Cl{\'\i}maco.
\newblock A review of interactive methods for multiobjective integer and
  mixed-integer programming.
\newblock \emph{European Journal of Operational Research}, 180\penalty0
  (1):\penalty0 99--115, 2007.

\bibitem[Angilella et~al.(2010)Angilella, Greco, and
  Matarazzo]{angilella2010non}
S.~Angilella, S.~Greco, and B.~Matarazzo.
\newblock Non-additive robust ordinal regression: A multiple criteria decision
  model based on the {C}hoquet integral.
\newblock \emph{European Journal of Operational Research}, 201\penalty0
  (1):\penalty0 277--288, 2010.

\bibitem[Berman et~al.(2010)Berman, Drezner, and Krass]{Berman:2010}
O.~Berman, Z.~Drezner, and D.~Krass.
\newblock Generalized coverage: new developments in covering location models.
\newblock \emph{Computer \& Operations Research}, 37\penalty0 (10):\penalty0
  1675--1687, 2010.

\bibitem[Bhattacharya and Bandyopadhyay(2010)]{bhattacharya2010solving}
R.~Bhattacharya and S.~Bandyopadhyay.
\newblock {Solving conflicting bi-objective facility location problem by {NSGA
  II} evolutionary algorithm}.
\newblock \emph{The International Journal of Advanced Manufacturing
  Technology}, 51\penalty0 (1-4):\penalty0 397--414, 2010.

\bibitem[Blanquero and Carrizosa(2002)]{Blanquero:2002}
R.~Blanquero and E.~Carrizosa.
\newblock A {DC} biobjective location model.
\newblock \emph{Journal of Global Optimization}, 23\penalty0 (2):\penalty0
  139--154, 2002.

\bibitem[Branke et~al.(2008)Branke, Deb, Miettinen, and
  S{\l}owi{\'n}ski]{BDMS_2008}
J.~Branke, K.~Deb, K.~Miettinen, and R.~S{\l}owi{\'n}ski, editors.
\newblock \emph{Multiobjective Optimization: Interactive and Evolutionary
  Approaches}, volume 5252 of \emph{LNCS}.
\newblock Springer, Berlin, 2008.

\bibitem[Branke et~al.(2015)Branke, Greco, S{\l}owi{\'n}ski, and
  Zielniewicz]{BGSZ-IEEE-TEC-2014}
J.~Branke, S.~Greco, R.~S{\l}owi{\'n}ski, and P.~Zielniewicz.
\newblock Learning {V}alue {F}unctions in {I}nteractive {E}volutionary
  {M}ultiobjective {O}ptimization.
\newblock \emph{IEEE Transactions on Evolutionary Computation}, 19\penalty0
  (1):\penalty0 88--102, 2015.

\bibitem[Branke et~al.(2016)Branke, Corrente, Greco, S{\l}owi{\'{n}}ski, and
  Zielniewicz]{NEMOIICh}
J.~Branke, S.~Corrente, S.~Greco, R.~S{\l}owi{\'{n}}ski, and P.~Zielniewicz.
\newblock Using {C}hoquet integral as preference model in interactive
  evolutionary multiobjective optimization.
\newblock \emph{European Journal of Operational Research}, 250:\penalty0
  884--901, 2016.

\bibitem[Branke et~al.(2017)Branke, Corrente, Greco, and Gutjahr]{bcggCOR}
J.~Branke, S.~Corrente, S.~Greco, and W.J. Gutjahr.
\newblock Efficient pairwise preference elicitation allowing for indifference.
\newblock \emph{Computers and Operations Research}, 88:\penalty0 175--186,
  2017.

\bibitem[Calik et~al.(2015)Calik, Labb{\'e}, and Yaman]{CalikLabbeYaman2015}
H.~Calik, M.~Labb{\'e}, and H.~Yaman.
\newblock {$p$-Center problems}.
\newblock In \emph{Location Science}, pages 79--92. Springer, 2015.

\bibitem[Carrizosa et~al.(2015)Carrizosa, Ushakov, and
  Vasilyev]{Carrizosa:2015}
E.~Carrizosa, A.~Ushakov, and I.~Vasilyev.
\newblock Threshold robustness in discrete facility location problems: a
  bi-objective approach.
\newblock \emph{Optimization Letters}, 9\penalty0 (7):\penalty0 1297--1314,
  2015.

\bibitem[Choquet(1953)]{choquet1953theory}
G.~Choquet.
\newblock Theory of capacities.
\newblock \emph{Annales de l'Institut Fourier}, 5\penalty0 (54):\penalty0
  131--295, 1953.

\bibitem[Church and ReVelle(1974)]{Church:1974}
R.L. Church and C.S. ReVelle.
\newblock The maximal covering location problem.
\newblock \emph{Papers in Regional Science}, 32\penalty0 (1):\penalty0
  101--118, 1974.

\bibitem[Ciomek et~al.(2017)Ciomek, Kadzi{\'n}ski, and
  Tervonen]{CiomekKadzinskiTervonen:2017}
K.~Ciomek, M.~Kadzi{\'n}ski, and T.~Tervonen.
\newblock Heuristics for prioritizing pair-wise elicitation questions with
  additive multi-attribute value models.
\newblock \emph{Omega}, 71:\penalty0 27--45, 2017.

\bibitem[Coello(2002)]{CoCo02}
C.A.~Coello Coello.
\newblock Theoretical and numerical constraint-handling techniques used with
  evolutionary algorithms: A survey of the state of the art.
\newblock \emph{Computer Methods in Applied Mechanics and Engineering},
  191\penalty0 (11-12):\penalty0 1245--1287, 2002.

\bibitem[Coutinho-Rodrigues et~al.(2012)Coutinho-Rodrigues, Tralh{\~a}o, and
  Al{\c{c}}ada-Almeida]{coutinho2012bi}
J.~Coutinho-Rodrigues, L.~Tralh{\~a}o, and L.~Al{\c{c}}ada-Almeida.
\newblock A bi-objective modeling approach applied to an urban semi-desirable
  facility location problem.
\newblock \emph{European Journal of Operational Research}, 223\penalty0
  (1):\penalty0 203--213, 2012.

\bibitem[Daskin(1995)]{Daskin:1995}
M.S. Daskin.
\newblock \emph{Network and discrete location: models, algorithms, and
  applications}.
\newblock Wiley, New York, USA, 1995.

\bibitem[Deb(2001)]{deb-book-01}
K.~Deb.
\newblock \emph{Multi-objective optimization using evolutionary algorithms}.
\newblock Chichester, UK: Wiley, 2001.

\bibitem[Deb et~al.(2002)Deb, Agrawal, Pratap, and Meyarivan]{debieee}
K.~Deb, S.~Agrawal, A.~Pratap, and T.~Meyarivan.
\newblock A fast and elitist multi-objective genetic algorithm: {NSGA-II}.
\newblock \emph{{IEEE} {T}ransactions on {E}volutionary {C}omputation},
  6\penalty0 (2):\penalty0 182--197, 2002.

\bibitem[Dias et~al.(2008)Dias, Captivo, and Cl{\'\i}maco]{dias2008memetic}
J.~Dias, M.E. Captivo, and J.~Cl{\'\i}maco.
\newblock A memetic algorithm for multi-objective dynamic location problems.
\newblock \emph{Journal of Global Optimization}, 42\penalty0 (2):\penalty0
  221--253, 2008.

\bibitem[Doerner et~al.(2009)Doerner, Gutjahr, and Nolz]{doerner2009multi}
K.F. Doerner, W.J. Gutjahr, and P.C. Nolz.
\newblock Multi-criteria location planning for public facilities in
  tsunami-prone coastal areas.
\newblock \emph{Or Spectrum}, 31\penalty0 (3):\penalty0 651--678, 2009.

\bibitem[Dom{\'\i}nguez-Mar{\'\i}n(2013)]{DominguezMarin2013}
P.~Dom{\'\i}nguez-Mar{\'\i}n.
\newblock \emph{{The Discrete Ordered Median Problem: Models and Solution
  Methods}}.
\newblock Springer Science \& Business Media, 2013.

\bibitem[Dom{\'\i}nguez-Mar{\'\i}n et~al.(2005)Dom{\'\i}nguez-Mar{\'\i}n,
  Nickel, Hansen, and Mladenovi{\'c}]{DominguezMarinEtAl2005}
P.~Dom{\'\i}nguez-Mar{\'\i}n, S.~Nickel, P.~Hansen, and N.~Mladenovi{\'c}.
\newblock {Heuristic procedures for solving the discrete ordered median
  problem}.
\newblock \emph{Annals of Operations Research}, 136\penalty0 (1):\penalty0
  145--173, 2005.

\bibitem[Drezner(2004)]{drezner2004location}
T.~Drezner.
\newblock Location of casualty collection points.
\newblock \emph{Environment and Planning C: Government and Policy}, 22\penalty0
  (6):\penalty0 899--912, 2004.

\bibitem[Drezner et~al.(2006)Drezner, Drezner, and Salhi]{drezner2006multi}
T.~Drezner, Z.~Drezner, and S.~Salhi.
\newblock A multi-objective heuristic approach for the casualty collection
  points location problem.
\newblock \emph{Journal of the Operational Research Society}, 57\penalty0
  (6):\penalty0 727--734, 2006.

\bibitem[Drezner and Hamacher(2001)]{Drezner:2004a}
Z.~Drezner and H.M. Hamacher.
\newblock \emph{Facility location: applications and theory}.
\newblock Springer, New York, USA, 2001.

\bibitem[Ehrgott and Gandibleux(2000)]{Ehrgott:2000}
M.~Ehrgott and X.~Gandibleux.
\newblock A survey and annotated bibliography of multiobjective combinatorial
  optimization.
\newblock \emph{OR Spectrum}, 22\penalty0 (4):\penalty0 425--460, 2000.

\bibitem[Ehrgott and Gandibleux(2008)]{Ehrgott:2008}
M.~Ehrgott and X.~Gandibleux.
\newblock Hybrid {M}etaheuristics for {M}ulti-objective {C}ombinatorial
  {O}ptimization.
\newblock \emph{Studies in Computational Intelligence}, 114:\penalty0 221--259,
  2008.

\bibitem[Eiben and Smith(2003)]{Eibe03a}
A.E. Eiben and J.E. Smith.
\newblock \emph{{Introduction to Evolutionary Computing}}.
\newblock Springer, 2003.

\bibitem[Eiselt and Laporte(1995)]{eiseltLaporte:1995}
H.A. Eiselt and G.~Laporte.
\newblock \emph{Objectives in location problems}.
\newblock Springer-Verlag, New York, USA, 1995.

\bibitem[Eiselt and Marianov(2011)]{knuth:art}
H.A. Eiselt and V.~Marianov.
\newblock \emph{Foundations of location analysis}.
\newblock International Series in Operations Research and Management Science.
  Springer, New York, USA, 2011.

\bibitem[Farahani et~al.(2010)Farahani, SteadieSeifi, and
  Asgari]{Farahanisurevey:2010}
R.~Z. Farahani, M.~SteadieSeifi, and N.~Asgari.
\newblock Multiple criteria facility location problems: A survey.
\newblock \emph{Operations Research}, 34\penalty0 (7):\penalty0 1689--1709,
  2010.

\bibitem[Fernandes et~al.(2014)Fernandes, Captivo, and
  Cl\'{\i}maco]{DSSCaptivo:2014}
S.~Fernandes, M.E. Captivo, and J.~Cl\'{\i}maco.
\newblock A {DSS} for bicriteria location problems.
\newblock \emph{Decision Support Systems}, 57:\penalty0 224--244, 2014.

\bibitem[Grabisch(1996)]{Grabisch1996}
M.~Grabisch.
\newblock The application of fuzzy integrals in multicriteria decision making.
\newblock \emph{European Journal of Operational Research}, 89\penalty0
  (3):\penalty0 445--456, 1996.

\bibitem[Grabisch(1997)]{grabisch1997k}
M.~Grabisch.
\newblock {$k$-order additive discrete fuzzy measures and their
  representation}.
\newblock \emph{Fuzzy sets and systems}, 92\penalty0 (2):\penalty0 167--189,
  1997.

\bibitem[Grabisch and Labreuche(2010)]{Grabisch2008}
M.~Grabisch and C.~Labreuche.
\newblock A decade of application of the {C}hoquet and {S}ugeno integrals in
  multi-criteria decision aid.
\newblock \emph{Annals of Operations Research}, 175\penalty0 (1):\penalty0
  247--290, 2010.

\bibitem[Greco et~al.(2016)Greco, Ehrgott, and
  Figueira]{GrecoEhrgottFigueira2016}
S.~Greco, M.~Ehrgott, and J.R. Figueira.
\newblock \emph{Multiple Criteria Decision Analysis: State of the Art Surveys}.
\newblock Springer, New York, 2016.

\bibitem[Hakimi(1964)]{Hakimi:1964}
S.L. Hakimi.
\newblock Optimum location of switching center and the absolute centers and
  medians of a graph.
\newblock \emph{Operations Research}, 12\penalty0 (3):\penalty0 450--459, 1964.

\bibitem[Hamacher et~al.(2002)Hamacher, Labbe, Nickel, and
  Skriver]{Hamacher:2002}
H.W. Hamacher, M.~Labbe, S.~Nickel, and A.J. Skriver.
\newblock Multicriteria semi-obnoxious network location problems ({MSNLP}) with
  sum and center objectives.
\newblock \emph{Annals of Operations Research}, 110\penalty0 (1-4):\penalty0
  33--53, 2002.

\bibitem[Harris et~al.(2011)Harris, Mumford, and Naim]{harris2011evolutionary}
I.~Harris, C.L. Mumford, and M.M. Naim.
\newblock {An evolutionary bi-objective approach to the capacitated facility
  location problem with cost and {CO2} emissions}.
\newblock In \emph{Proceedings of the 13th annual conference on Genetic and
  evolutionary computation}, pages 697--704. ACM, 2011.

\bibitem[Heyns and van Vuuren(2015)]{heyns2015multi}
AM~Heyns and JH~van Vuuren.
\newblock {Multi-objective optimisation of discrete {GIS}-based facility
  location problems}.
\newblock \emph{Optimization and Engineering}, 2015.

\bibitem[Hollander et~al.(2013)Hollander, Wolfe, and
  Chicken]{HollanderWolfeChicken2013}
M.~Hollander, D.A. Wolfe, and E.~Chicken.
\newblock \emph{Nonparametric statistical methods}, volume 751.
\newblock John Wiley \& Sons, 2013.

\bibitem[Jacquet-Lagreze and Siskos(1982)]{jacquet1982assessing}
E.~Jacquet-Lagreze and Y.~Siskos.
\newblock Assessing a set of additive utility functions for multicriteria
  decision-making, the {UTA} method.
\newblock \emph{European Journal of Operational Research}, 10\penalty0
  (2):\penalty0 151--164, 1982.

\bibitem[Kalcsics et~al.(2014)Kalcsics, Nickel, Pozo, Puerto, and
  Rodr\'{i}guez-Ch\'{i}a]{Pozo:2014}
J.~Kalcsics, S.~Nickel, M.A. Pozo, J.~Puerto, and A.M. Rodr\'{i}guez-Ch\'{i}a.
\newblock The multicriteria $p$-facility median location problem on networks.
\newblock \emph{European Journal of Operational Research}, 235\penalty0
  (3):\penalty0 484--493, 2014.

\bibitem[Karasakal and Nadirler(2008)]{karasakal2008interactive}
E.~Karasakal and D.~Nadirler.
\newblock An interactive solution approach for a bi-objective semi-desirable
  location problem.
\newblock \emph{Journal of Global Optimization}, 42\penalty0 (2):\penalty0
  177--199, 2008.

\bibitem[Kariv and Hakimi(1969)]{Kariv:1969}
O.~Kariv and S.L. Hakimi.
\newblock {An algorithmic approach to network location problems. {I}: the
  $p$-medians}.
\newblock \emph{SIAM Journal of Applied Mathematics}, 37\penalty0 (3):\penalty0
  539--560, 1969.

\bibitem[Keeney and Raiffa(1976)]{Keeney76}
{R.L}. Keeney and H.~Raiffa.
\newblock \emph{Decisions with multiple objectives: Preferences and value
  tradeoffs}.
\newblock J. Wiley, New York, 1976.

\bibitem[Knowles et~al.(2001)Knowles, Watson, and
  Corne]{KnowlesWatsonCorne2001}
J.D. Knowles, R.A. Watson, and D.W. Corne.
\newblock {Reducing local optima in single-objective problems by
  multi-objectivization}.
\newblock In \emph{International conference on evolutionary Multi-Criterion
  Optimization}, pages 269--283. Springer, 2001.

\bibitem[Krarup and Pruzan(1983)]{krarup1983simple}
J.~Krarup and P.M. Pruzan.
\newblock The simple plant location problem: survey and synthesis.
\newblock \emph{European Journal of Operational Research}, 12\penalty0
  (1):\penalty0 36--81, 1983.

\bibitem[Laporte et~al.(2015)Laporte, Nickel, and da~Gama]{LocationScience}
G.~Laporte, S.~Nickel, and F.S. da~Gama.
\newblock \emph{Location science}.
\newblock Springer, Berlin, 2015.

\bibitem[Marsh and Schilling(1994)]{Marsh:1994}
M.T. Marsh and D.A. Schilling.
\newblock {Equity measurement in facility location analysis: A review and
  framework}.
\newblock \emph{European Journal of Operational Research}, 74\penalty0
  (1):\penalty0 1--7, 1994.

\bibitem[Miettinen et~al.(2008)Miettinen, Ruiz, and
  Wierzbicki]{miettinen2008introduction}
K.~Miettinen, F.~Ruiz, and A.P. Wierzbicki.
\newblock Introduction to multiobjective optimization: interactive approaches.
\newblock In J.~Branke, K.~Deb, R.~S{\l}owi\'{n}ski, and K.~Miettinen, editors,
  \emph{Multiobjective optimization}, pages 27--57. Berlin: Springer, 2008.

\bibitem[Mladenovic et~al.(2007)Mladenovic, Brimberg, Hansen, and
  Moreno-Perez]{Mlad:2007}
N.~Mladenovic, J.~Brimberg, P.~Hansen, and J.A. Moreno-Perez.
\newblock {The $p$-median problem: A survey of metaheuristic approaches}.
\newblock \emph{European Journal of Operational Research}, 179\penalty0
  (3):\penalty0 927--939, 2007.

\bibitem[Nelder and Mead(1965)]{nelder1965simplex}
J.A. Nelder and R.~Mead.
\newblock A simplex method for function minimization.
\newblock \emph{The computer journal}, 7\penalty0 (4):\penalty0 308--313, 1965.

\bibitem[Nickel(2001)]{Nickel2001}
S.~Nickel.
\newblock {Dicrete Ordered Weber problems}.
\newblock In R.~Fleischmann, U.~Lasch, U.~Derigs, W.~Domschke, and U.~Rieder,
  editors, \emph{{Operations Research Proceedings 2000}}, pages 71--76.
  Springer, 2001.

\bibitem[Nijkamp and Spronk(1981)]{nijkamp1981interactive}
P.~Nijkamp and J.~Spronk.
\newblock Interactive multidimensional programming models for locational
  decisions.
\newblock \emph{European Journal of Operational Research}, 6\penalty0
  (2):\penalty0 220--223, 1981.

\bibitem[Ohsawa and Tamura(2003)]{Ohsawa:2003}
Y.~Ohsawa and K.~Tamura.
\newblock Efficient location for a semi-obnoxious facility.
\newblock \emph{Annals of Operations Research}, 123\penalty0 (1-4):\penalty0
  173--188, 2003.

\bibitem[Ohsawa et~al.(2008)Ohsawa, Ozaki, and Plastria]{Ohsawa:2007}
Y.~Ohsawa, N.~Ozaki, and F.~Plastria.
\newblock {Equity-effciency bicriteria location with squared Euclidean
  distances}.
\newblock \emph{Operations Research}, 56\penalty0 (1):\penalty0 79--87, 2008.

\bibitem[Owen and Daskin(1998)]{Daskin}
S.H. Owen and M.S. Daskin.
\newblock Strategic facility location: a review.
\newblock \emph{European Journal of Operational Research}, 111\penalty0
  (3):\penalty0 423--447, 1998.

\bibitem[Rahmati et~al.(2014)Rahmati, Ahmadi, Sharifi, and
  Chambari]{rahmati2014multi}
S.H.A. Rahmati, A.~Ahmadi, M.~Sharifi, and A.~Chambari.
\newblock A multi-objective model for facility location--allocation problem
  with immobile servers within queuing framework.
\newblock \emph{Computers \& Industrial Engineering}, 74:\penalty0 1--10, 2014.

\bibitem[Rakas et~al.(2004)Rakas, Teodorovi\'{c}, and Kim]{Rakas:2004}
J.~Rakas, D.~Teodorovi\'{c}, and T.~Kim.
\newblock Multi-objective modeling for determining location of undesirable
  facilities.
\newblock \emph{Transportation Research Part D: Transport and Environment},
  9\penalty0 (2):\penalty0 125--138, 2004.

\bibitem[ReVelle and Swaim(1970)]{revelleswaim:1970}
C.S. ReVelle and R.W. Swaim.
\newblock Central facilities location.
\newblock \emph{Geographical Analysis}, 2\penalty0 (1):\penalty0 30--42, 1970.

\bibitem[Rota(1964)]{Rota}
G.C. Rota.
\newblock On the foundations of combinatorial theory {I.} {T}heory of
  {M}\"{o}bius functions.
\newblock \emph{Wahrscheinlichkeitstheorie und Verwandte Gebiete}, 2:\penalty0
  340--368, 1964.

\bibitem[Roy(1987)]{roy1987meaning}
B.~Roy.
\newblock Meaning and validity of interactive procedures as tools for decision
  making.
\newblock \emph{European Journal of Operational Research}, 31\penalty0
  (3):\penalty0 297--303, 1987.

\bibitem[Roy(1993)]{roy1993decision}
B.~Roy.
\newblock Decision science or decision-aid science?
\newblock \emph{European Journal of Operational Research}, 66\penalty0
  (2):\penalty0 184--203, 1993.

\bibitem[Roy(2016)]{Roy_book_greco_2016}
B.~Roy.
\newblock Paradigm and {C}hallenges.
\newblock In S.~Greco, J.R. Figueira, and M.~Ehrgott, editors, \emph{{Multiple
  Criteria Decision Analysis: State of the Art Surveys}}, pages 19--39.
  Springer, New York, 2016.

\bibitem[Saaty(1977)]{Saaty1977}
T.~Saaty.
\newblock A scaling method for priorities in hierarchical structures.
\newblock \emph{Journal of Mathematical Psychology}, 15\penalty0 (3):\penalty0
  234--281, 1977.

\bibitem[Saaty(1980)]{Saaty1980}
T.~Saaty.
\newblock \emph{The {A}nalytic {H}ierarchy {P}rocess}.
\newblock New York, McGraw-Hill, 1980.

\bibitem[Schnepper et~al.(2019)Schnepper, Klamroth, Stiglmayr, and
  Puerto]{SchnepperEtAl2019}
T.~Schnepper, K.~Klamroth, M.~Stiglmayr, and J.~Puerto.
\newblock {Exact algorithms for handling outliers in center location problems
  on networks using $k$-max functions}.
\newblock \emph{European Journal of Operational Research}, 273\penalty0
  (2):\penalty0 441--451, 2019.

\bibitem[Segura et~al.(2016)Segura, Coello~Coello, Miranda, and
  Le{\'o}n]{SeguraEtAl2016}
C.~Segura, C.A. Coello~Coello, G.~Miranda, and C.~Le{\'o}n.
\newblock {Using multi-objective evolutionary algorithms for single-objective
  constrained and unconstrained optimization}.
\newblock \emph{Annals of Operations Research}, 240\penalty0 (1):\penalty0
  217--250, 2016.

\bibitem[Serafini(1987)]{Serafini1987}
P.~Serafini.
\newblock Some considerations about computational complexity for multi
  objective combinatorial problems.
\newblock In J.~Jahn and W.~Krabs, editors, \emph{Recent Advances and
  Historical Development of Vector Optimization}, pages 222--232. Springer,
  1987.

\bibitem[Shafer(1976)]{Shafer}
G.~Shafer.
\newblock \emph{A Mathematical Theory of Evidence}.
\newblock Princeton University Press, 1976.

\bibitem[Shankar et~al.(2013)Shankar, Basavarajappa, Chen, and
  Kadadevaramath]{shankar2013location}
B.L. Shankar, S.~Basavarajappa, J.C.H. Chen, and R.S. Kadadevaramath.
\newblock {Location and allocation decisions for multi-echelon supply chain
  network--A multi-objective evolutionary approach}.
\newblock \emph{Expert Systems with Applications}, 40\penalty0 (2):\penalty0
  551--562, 2013.

\bibitem[Tansel et~al.(1982)Tansel, Francis, and Lowe]{Tansel:1982}
B.C. Tansel, R.L. Francis, and T.J. Lowe.
\newblock A biobjective multifacility minimax location problem on a tree
  network.
\newblock \emph{Transportation Science}, 16\penalty0 (4):\penalty0 407--429,
  1982.

\bibitem[Teghem et~al.(2000)Teghem, Tuyttens, and
  Ulungu]{teghem2000interactive}
J.~Teghem, D.~Tuyttens, and E.L. Ulungu.
\newblock An interactive heuristic method for multi-objective combinatorial
  optimization.
\newblock \emph{Computers \& Operations Research}, 27\penalty0 (7):\penalty0
  621--634, 2000.

\bibitem[Tomczyk and Kadzinski(2019)]{TomczykKadzinski2019}
M.K. Tomczyk and M.~Kadzinski.
\newblock {EMOSOR: E}volutionary multiple objective optimization guided by
  interactive stochastic ordinal regression.
\newblock \emph{Computers \& Operations Research}, 108:\penalty0 134 -- 154,
  2019.

\bibitem[Vaidya and Kumar(2006)]{VaidyaKumar2006}
O.S. Vaidya and S.l Kumar.
\newblock Analytic hierarchy process: {A}n overview of applications.
\newblock \emph{European Journal of operational research}, 169\penalty0
  (1):\penalty0 1--29, 2006.

\bibitem[Van~Veldhuizen(1999)]{VanVeldhuizenPhD}
D.~Van~Veldhuizen.
\newblock Multiobjective evolutionary algorithms: Classifications, analysis and
  new innovations.
\newblock Ph{D} thesis, {A}ir {F}orce {I}nstitute of {T}echnology, {F}aculty of
  the {G}raduate school of {E}ngineering, 1999.

\bibitem[Villegas et~al.(2006)Villegas, Palacios, and
  Medaglia]{villegas2006solution}
J.G. Villegas, F.~Palacios, and A.L. Medaglia.
\newblock Solution methods for the bi-objective (cost-coverage) unconstrained
  facility location problem with an illustrative example.
\newblock \emph{Annals of Operations Research}, 147\penalty0 (1):\penalty0
  109--141, 2006.

\bibitem[Wakker(1989)]{wakker1989additive}
P.P. Wakker.
\newblock \emph{Additive representations of preferences: A new foundation of
  decision analysis}.
\newblock Springer, 1989.

\bibitem[Yapicioglu et~al.(2007)Yapicioglu, Smith, and
  Dozier]{yapicioglu2007solving}
H.~Yapicioglu, A.E. Smith, and G.~Dozier.
\newblock Solving the semi-desirable facility location problem using
  bi-objective particle swarm.
\newblock \emph{European Journal of Operational Research}, 177\penalty0
  (2):\penalty0 733--749, 2007.

\bibitem[Zhou et~al.(2011)Zhou, Qu, Li, Zhao, Suganthan, and
  Zhang]{zhou2011multiobjective}
A.~Zhou, B.Y. Qu, H.~Li, S.Z. Zhao, P.N. Suganthan, and Q.~Zhang.
\newblock Multiobjective evolutionary algorithms: {A} survey of the state of
  the art.
\newblock \emph{Swarm and Evolutionary Computation}, 1\penalty0 (1):\penalty0
  32--49, 2011.

\bibitem[Zitzler et~al.(2002)Zitzler, Laumanns, and Thiele]{zlt2002a}
E.~Zitzler, M.~Laumanns, and L.~Thiele.
\newblock {SPEA2: Improving the Strength Pareto Evolutionary Algorithm} for
  multiobjective optimization.
\newblock In K.C. Giannakoglou et~al., editors, \emph{Evolutionary Methods for
  Design, Optimisation and Control with Application to Industrial Problems
  (EUROGEN 2001)}, pages 95--100. International Center for Numerical Methods in
  Engineering (CIMNE), 2002.

\end{thebibliography}

\end{document}